\documentclass[12pt,leqno]{article}
\usepackage{amsmath,amssymb, amscd}
\newcommand{\be}{\begin{equation}}
\newcommand{\ee}{\end{equation}}
\newcommand{\bea}{\begin{eqnarray}}
\newcommand{\eea}{\end{eqnarray}}
\newcommand{\bean}{\begin{eqnarray*}}
\newcommand{\eean}{\end{eqnarray*}}

\newcommand{\Q}{\mathbb{Q}}

\newcommand{\N}{\mathbb{N}}
\newcommand{\zz}{\mathbb{Z}}

\newcommand{\ga}{\alpha}

\newcommand{\pp}{\mathbb{P}}

\newenvironment{remark}[1][Remark]{\begin{trivlist}
       \item[\hskip \labelsep {\bfseries #1}]}{\end{trivlist}}

\newenvironment{caja}[1]{\begin{trivlist}
       \item[\hskip \labelsep {\bfseries #1}]}{\end{trivlist}}

\begin{document}

\title{\bf On the Lefschetz Standard Conjecture}
\author{Jos\'e J. Ram\'on Mar\'i}
\date{ }

\maketitle
\newtheorem{theorem}{Theorem}[section]
\newtheorem{prop}[theorem]{Proposition}
\newtheorem{corollary}[theorem]{Corollary}
\newtheorem{lemma}[theorem]{Lemma}
\newtheorem{obs}[theorem]{Remark}

\abstract{The subject of the present paper is Grothendieck's
Lefschetz standard conjecture $B(X)$.  Our main result is that, if
$X$ is a projective smooth variety of dimension $n$ and the
conjecture $B({\cal Y})$ holds for the generic fibre ${\cal Y}$
(of dimension $n-1$ over the field $k(t)$) of a suitable Lefschetz
fibration of $X$, then the operator $\Lambda_X-p^{n+1}_X$ is
algebraic.  If in addition $p^{n+1}_X$ is algebraic, then $B(X)$
is settled. Along the way we establish the algebraicity of the K\"
unneth projectors $\pi^i_X$ for $i\neq n-1, n, n+1$ under the
above hypotheses.}

\vspace*{1cm}

\section{Introduction}

  All varieties involved are assumed to be smooth and projective, unless otherwise stated.  The notations on correspondences that we adopt are those of Kleiman ~\cite{Kleiman} 1.3, Jannsen ~\cite{Jannsen}, Scholl ~\cite{Scholl}.  We fix a prime $\ell\neq \mbox{char }k$, and make the harmless assumption that our field $k$, of arbitrary characteristic, contains the $\ell^N$-th roots of unity for all $N$; then we go on with the `heresy'  ~\cite{G} $\zz_\ell(1) \approx \zz_\ell$.

  Let $X$ be a smooth projective variety of dimension $n$ over a field $k$; we now fix a very ample line bundle ${\cal L}$, giving an immersion in
  $\pp^N$, which we replace if necessary by a tensor power ${\cal L}^{\otimes m}$ to obtain condition {\bf (A)} of Section \ref{lapiceslefschetz}.  Let $Y$ be a smooth hyperplane section; we write $\xi_X:=[Y]\in H^2(X)(1)$.  Let $L_X$ ($L$ when not misleading) be the Lefschetz operator $L_Xx=[Y]\wedge x$, where $\wedge$ denotes the cup-product in $H^*(X)$.  $A^*(X)$ will denote the graded ring of algebraic cycles modulo homological equivalence with coefficients over $\Q$,
   and $A^{n+*}(X\times X)$ will denote the ring of homological correspondences (with coefficients over $\Q$), $\circ$ being the product considered.
     Note that the degree of a correspondence $u\in CH^{\mathrm{dim }\; X+r}(X\times X')$ is $r$ as usual ~\cite{F}, and the  cohomological degree of $u$, \textit{i.e.}  the degree of $u$ as an operator in cohomology  $H^*(X)\rightarrow H^*(Y)$ will be
  $2r$. Given a subspace $V$ of $H^*(X)$, we denote by $e_V$ the
  orthogonal projection onto $V$.
   Following Kleiman ~\cite{Kleiman}, we denote the trace (or orientation) map by $\langle \rangle:H^*(X)\rightarrow \Q_\ell$ and the Poincar\'e duality

    pairing by $\langle , \rangle:H^i(X)\otimes H^{2n-i}(X)\rightarrow \Q_\ell$.

  The Hard Lefschetz Theorem ~\cite{DeWeII} states that the maps
  $$L^{n-i}:H^i(X)\rightarrow H^{2n-i}(X)$$ are isomorphisms (henceforth called \textit{Lefschetz isomorphisms}).  One
  can then define the primitive subspaces $P^i(X)=\mbox{Ker
  }L^{n-i+1}\cap H^i(X)$, and one has a Lefschetz decomposition of
  $H^*(X)$: $H^i(X)=\oplus L^jP^{i-2j}(X)$.  Let Let $x=\sum L^j x_{i-2j}$ be the Lefschetz decomposition of $x \in H^i(X).$  Denote $i_1=max\{i-n,1\}$. We define the following
operators of degree $-2$:
$$\Lambda x= \sum_{j\geq i_1} L^{j-1}x_{i-2j},$$
$$^c\Lambda x=\sum_{j\geq i_1} j(n-i+j+1) L^{j-1} x_{i-2j}.$$
 We denote the K\"unneth projectors $H^*(X)\twoheadrightarrow H^i(X) \hookrightarrow H^*(X)$ by $\pi_X^i=\pi^i.$ We define the operator of degree $0$
  $$H=H_X=\sum_{i=0}^{2n} (n-i)\pi_X^i.$$
    The following operators are also essential: for $x=\sum L^j x_{i-2j}\in H^i(X)$, $p^kx=\delta_{i,k} x_k$, when $i\leq n$,
     and $p^kx=\delta_{i,k}x_{2n-k}$ for $k>n$; it is clear that $p^i$ is a projector for $i \leq n$.
Whenever we have polarised varieties $X_i$, we will consider the
induced polarisation on $X_1\times X_2$, and so $L_{X_1\times
X_2}=L_{X_1}\otimes 1+ 1\otimes L_{X_2}.$  We will do likewise
when we have an inclusion; for instance, let $\iota:Y\subset X$
denote an inclusion of a smooth hyperplane section.  Then
$\xi_Y=\iota^*\xi_X$ and $L_X=\iota_*\iota^*, L_Y=\iota^*\iota_*$.
We denote the space of vanishing cycles by $V(Y)=\mbox{Ker
}\iota_*|H^{n-1}(Y) \subset H^{n-1}(Y)$, with $Y$ as above.

  We recall the following result:
\begin{prop}\label{sl2}(Kleiman ~\cite{Kleiman} 1.4.6, ~\cite{Andre}) The operators $ ^c\Lambda, L,H$ are an ${\mathfrak sl}_2$- triple; in other words, the following identities hold:
$$[ ^c\Lambda,L]=H, \;[H,L]=-2L,\; [H, ^c\Lambda]=2 ^c\Lambda.$$\end{prop}

The following conjecture was stated by Grothendieck, and is one of
his \emph{standard conjectures} ~\cite{G} ~\cite{Kleiman}:
\begin{caja}{B(X):} The operator $\Lambda$ is induced by an algebraic cycle; equivalently (~\cite{Kleiman} Prop. 2.3), all the operators in the ${\mathfrak sl}_2$-triple $(^c\Lambda,L,H)$ are algebraic.\end{caja}
  The conjecture $B(X)$ is known for curves, surfaces, generalised flag varieties, abelian
  varieties and is stable under products and smooth hyperplane
  sections ~\cite{Kleiman}.  We will therefore assume that $n\geq 3.$
 For a discussion on this form of the conjecture --regarding the field of definition-- see \ref{finalcomments}.  Another standard conjecture
 of Grothendieck,  weaker than $B(X)$,  (Kleiman ~\cite{Kleiman} 2.4, ~\cite{G}) regards the
algebraicity of the K\"unneth projectors (again, we refer to
\ref{finalcomments}):
\begin{caja}{C(X):} The K\"unneth projectors $\pi^i$ are algebraic for all $i=0,\ldots, 2n$.\end{caja}

  The main result of this paper states as follows.
  \begin{caja}{Main Theorem.} Let $X$ be smooth projective of dimension $n\geq 3$.  Assume the
  conjecture $B({\cal Y})$ for the general fibre ${\cal Y}$ of a Lefschetz pencil of $X$ satisfying condition {\bf (A)} (see Section \ref{lapiceslefschetz}).
  Then the operator $\Lambda_X-p^{n+1}_X$ is algebraic.
\end{caja}
The following result on $C(X)$ is proven in Proposition
\ref{bandc}, although it may be viewed as a corollary of the Main
Theorem (see subsection \ref{subsectmainth}):
\begin{caja}{Partial result on C(X).} Assume $B({\cal Y})$ for ${\cal Y}$ as above.  Then the K\"unneth projectors $\pi^i_X$ are
  algebraic for all $i\neq n-1, n, n+1$.\end{caja}

  The departure point of our proof is essentially the algebraic cycle $\Lambda_{\cal
  Y}$ on the generic fibre ${\cal Y}/k(t)$ of a Lefschetz fibration of
  $X$ satisfying condition {\bf (A)} of Section \ref{lapiceslefschetz} (Katz ~\cite{SGA7} XVIII.5.3),
  $\rho:\tilde{X}\rightarrow\pp^1.$  In our proof, we pay special
  attention to
  the correspondences supported on
  $\tilde{X}\times_{\pp^1}\tilde{X}$ (see \ref{relcorresps}),
  which turn out to preserve the Leray filtration of $\rho$,
   as will be seen in Proposition \ref{fil}.  Whenever we have an algebraic class
  $u$ in $A^{n-1+*}({\cal Y}\times {\cal Y})$, a \textbf{lifting} (or \textbf{extension}) of $u$ will denote a class supported on $\tilde{X}\times_{\pp^1}\tilde{X}$
   which yields $u$ after restriction to the generic fibre of $\pp^1$.  Our
   proof requires that the assumptions be over $k(t)$.

    The consequences of establishing the full conjecture $B(X)$
    for general $X$ would be
    remarkable.  Not only would this yield a satisfactory category of
    pure motives in characteristic zero, but as shown by  Y. Andr\'e
    ~\cite{Andre} it would imply the Variational Hodge Conjecture,
 hence the Hodge conjecture for arbitrary products of the form $A\times X_1\times \cdots \times
 X_m,$ where $A$ is an abelian variety and $X_i$ are K3 surfaces.
\begin{caja}{Acknowledgements. } My enormous debt to the work of A. Grothendieck is apparent.   This paper would never have been written were it
not for the unwavering support and reassurance I received from
Mar\'ia M. Morales, my \textit{madrina}.  I am indebted to the
following people for their continuous encouragement and support:
Brendan Guilfoyle, Luke `Dr Socks' Murray, Esteban Roig Davison,
Virgolina Lopes Holme.  The interest and support that Kevin
Hutchinson and Gary McGuire (UCD) showed is gratefully
acknowledged. Finally, I heartily thank David Simms for providing
with full access to the TCD Mathematics Library.
\end{caja}

\section{General results}
\label{seccgeneral}
 The results in this section need no more background than Kleiman ~\cite{Kleiman}.
  For the sake of completeness, we include the proof of the following lemma.
\begin{lemma}\label{adjuntoh} Let $u \in A^{n+r}(X\times X)$ be a correspondence of degree $r$ on $X$.  The following identity holds:
$$[H,u]=-2r\cdotp u.$$\end{lemma}
{\bf Proof:} One has $u\pi^i=\pi^{i+2r}u$, hence
$$uH=\sum (n-i)u\pi^i=\sum (n-i)\pi^{i+2r}u=Hu+2r\cdotp u.$$
 Isolating yields $[H,u]=-2r\cdotp u$ as desired. $\blacksquare$

\begin{lemma}\label{cgen}Let $f:X'\rightarrow X$ be a generically finite, surjective morphism of smooth projective varieties.  Assume that $C(X')$ holds; then $C(X)$ holds.\end{lemma}
 The lemma follows readily from the identity $\pi^i_X=\frac{1}{deg(f)}f_* \pi^i_{X'} f^*$. $\blacksquare$

\begin{lemma}\label{primis}  If $j_1+2i_1+j_2+2i_2=2n$, then the pieces $L^{i_1}P^{j_1}(X)$ and $L^{i_2}P^{j_2}(X)$ are orthogonal for $j_1 \neq j_2$.  Let $0\leq i\leq n$; then the operator $p^i$ is a projector, and $p^{2n-i}$ is a symmetric operator
characterised by $p^{2n-i}:L^{n-i}P^i(X)\rightarrow P^i(X)$ is
given by the inverse of the Lefschetz isomorphism $L^{n-i}$, and
$p^{2n-i}L^jP^k(X)=0$ if $(j,k)\neq (n-i,i).$  The following
identities hold:
$$p^{2n-i}L^{n-i}=p^i, \qquad L^{n-i}p^{2n-i}=\,^tp^i.$$\end{lemma}
{\bf Proof:} The first assertion implies the rest of the lemma.
Suppose $j_2>j_1$; then
$$L^{i_2}P^{j_2}(X)\wedge L^{i_1}P^{j_1}(X)=L^{i_1+i_2}P^{j_1}(X)\wedge P^{j_2}(X)=0;$$
indeed, consider a smooth linear section $\kappa:W \hookrightarrow
X$ of codimension $(i_1+i_2)$.  Then
$\kappa_*\kappa^*(P^{j_1}(X)\wedge P^{j_2}(X))=0,$ since
$j_2>\frac{j_1+j_2}{2}=\mbox{dim }W$, hence $\kappa^*P^{j_2}(X)=0$
by ~\cite{Kleiman} 1.4.7, ~\cite{SGA7} Exp. XVIII (5.2.4).
  This proves the assertion. $\blacksquare$

\begin{lemma} The operators $L, \Lambda$ and $^c\Lambda$ are symmetric.
\end{lemma}
{\bf Proof:} The first two are well-known ~\cite{Kleiman}. We
prove $^c\Lambda=\;^t(^c\Lambda)$.  By Lemma \ref{primis} one need
only check the following. Let $x=L^jx_{i-2j}\in L^jP^{i-2j}(X),
y=L^{n-i+j+1}y_{i-2j}\in L^{n-i+j+1}P^{i-2j}(X)$; the following
equality holds:
$$\langle  ^c\Lambda L^jx_{i-2j}, L^{n-i+j+1}y_{i-2j}\rangle=\langle L^jx_{i-2j}, ^c\Lambda L^{n-i+j+1}y_{i-2j} \rangle.$$
  Indeed, let $c=\langle L^{n-i} x_{i-2j}, y_{i-2j}\rangle$.  Then   $\mbox{l.h.s.}=j(n-i+j+1) c$ and
$\mbox{r.h.s.}=(n-i+j+1)j c,$ since $j=n-(2n-i+2)+(n-i+j+1)+1$. The Lemma is thus established. $\blacksquare$

  \begin{prop}\label{basiclambda}\begin{enumerate}
    \item The following non-commutative rings of operators are
    equal: $$\Q\langle L, \Lambda\rangle=\Q\langle L, \,^c\Lambda \rangle=\Q\langle L, p^n, \cdots , p^{2n}\rangle.$$
    \item $B(X)$ holds if and only if, for all $i<n$, the inverse
    $$\theta^i:H^{2n-i}(X)\rightarrow H^{i}(X)$$ to the Lefschetz isomorphism $L^{n-i}:H^i(X)\overset{\sim}{\rightarrow} H^{2n-i}(X)$
    is induced by an algebraic correspondence for $i<n$.
  \end{enumerate}
  \end{prop}
The first assertion follows from Kleiman ~\cite{Kleiman} Prop.
1.4.4.  The second is proved in \textit{op. cit.}, 2.3.
$\blacksquare$

  The morphisms $\iota^*, \iota_*$ are well-behaved with respect to the
  Lefschetz decompositions of $X, Y$ (see Kleiman ~\cite{Kleiman} Prop. 1.4.7).  The following two lemmas
  relate the operators $\Lambda_X, \Lambda_Y$.
\begin{lemma}\label{prinduccion} The following identity holds:
\be \iota^*\Lambda_X=\Lambda_Y\iota^*+\sum_{j=n+1}^{2n-2}
\iota^*L^{j-n-1} p^j_X. \ee
\end{lemma}
{\bf Proof:} Here we use ~\cite{Kleiman} Prop. 1.4.7 constantly.
The operators $\iota^*\Lambda_X$ and $\Lambda_Y\iota^*$ agree on
$H^i(X)$ for $i\leq n$.  It is easy to see that
$\iota^*\Lambda_X-\Lambda_Y\iota^*=0$ on $L^{i-n+1}H^{2n-i-2}(X)$.
The equality
$\iota^*\Lambda_X-\Lambda_Y\iota^*=\sum_{i=n+1}^{2n-2}\iota^*L^{j-n-1}
p^j_X$ thus holds for every piece $L^rP^s(X)$ of $H^*(X)$, and the
proof is now complete.
$\blacksquare$

  The following is but a rephrasing of Kleiman ~\cite{Kleiman} Prop. 2.12.
\begin{lemma}\label{iotalambda}  Assuming $B(Y)$, $B(X)$ is equivalent to $\iota^*\Lambda_X$ being algebraic.\end{lemma}
{\bf Proof:} We reproduce the proof in ~\cite{Kleiman}.  First
note that $ ^t(\iota^*\Lambda_X)=\Lambda_X\iota_*$.  Assume that
$\iota^*\Lambda_X$ is algebraic; then $\Lambda_Y=\iota^*\Lambda_X
\,^t(\iota^*\Lambda_X)=\iota^*\Lambda_X^2\iota_*$ is algebraic
and, for every $i<n$, the map

\be
\begin{CD}
\theta^i:H^{2n-i}(X)@>\iota^*\Lambda_X>>
H^{2n-2-i}(Y)@>\Lambda_Y^{n-1-i}>> H^i(Y)@>\Lambda_X\iota_*>>
H^i(X)
\end{CD}
\ee

is an algebraic inverse to $L^{n-i}$.  This establishes $B(X)$ by
 Proposition \ref{basiclambda}. $\blacksquare$
\vskip .5cm
   The following will be used in Section \ref{proyrelativos}.
  \begin{lemma}\label{cx} The conjecture $C(X)$ holds if and only if the semisimple operator $H$ is algebraic.\end{lemma}
{\bf Proof:} The result follows readily from the identities
$$[\Delta_X]=id_{H^*(X)}=\sum \pi^i \mbox{  and  } H^r=\sum
(n-i)^r\pi^i \mbox{  for  } r\in \N.$$ $\blacksquare$
\section{The cohomology of Lefschetz pencils}
    For the basic results and the tone of this section we follow Katz ~\cite{SGA7} Exp. XVIII; we assume $k$ to be algebraically closed.  For $X$ an $n$-dimensional variety and $X \hookrightarrow
  \pp^N$ a suitable embedding, there exists a line $L\subset
  (\pp^N)^\vee$ cutting the dual variety $X^\vee$ of $X\subset \pp^N$
  transversally; $L$ is then called a \textbf{Lefschetz pencil.}
   A basic property of $L$ is that, for every hyperplane $t \in L,$ $X_t=X\cap H_t$ is either smooth or has a unique singular point which is an
  ordinary double point.  The base locus of $L$ in $X$ will be
  denoted by $\Delta$, and for any two $t_1\neq t_2\in L$ one has
  a transversal intersection $X_{t_1}\cap X_{t_2}=\Delta$.  Thus
  $\Delta$ is smooth of dimension $n-2$: for any smooth member $Y=X_t$ as above, we will denote
  the canonical inclusion by $h:\Delta \hookrightarrow Y$.  If $\tilde{X}$ denotes
  the blowing-up of $X$ centred at $\Delta$, projection induces a
  map $X-\Delta\rightarrow \pp^1\cong L$ which induces a fibration
  (henceforth called a \textbf{Lefschetz fibration}, or \textbf{Lefschetz
  pencil} by \emph{abus de langage}):
  $$\rho:\tilde{X}\longrightarrow \pp^1.$$
  We denote by $f$ the blowing-up map $f:\tilde{X}\rightarrow X$.
  The full blow-up diagram will be denoted as follows:
  \be
\begin{CD}
\tilde{\Delta}  @>i>>   \tilde{X} \\
       @VVgV                @VVfV \\
        \Delta  @>j>>       X,
\end{CD}
\ee  where $j$ is the canonical inclusion $\Delta \subset X$.
 $\tilde{\Delta}$ is then the exceptional divisor and, $\Delta$
being a complete intersection, $g$ is a trivial projective bundle;
$j$ denotes the inclusion $\tilde{\Delta} \subset \tilde{X}.$  We
describe the cohomology of $\tilde{X}$ in the following
proposition.
\begin{prop} \label{katzuno}(Katz ~\cite{SGA7} Exp. XVIII Prop. 4.2) Notations and assumptions being as above.  Then: \begin{enumerate}
 \item[(i)] the following homomorphisms are mutual inverses:
 $$H^\bullet(\tilde{X})\overset{f_*\oplus g_*i^*}{\longrightarrow}
 H^\bullet(X)\oplus H^{\bullet -2}(\Delta)(-1)$$
 and
$$ H^\bullet(X)\oplus H^{\bullet
 -2}(\Delta)(-1)\overset{f^*+i_*g^*}{\longrightarrow}H^\bullet(\tilde{X}).$$
 \item[(ii)] Transport of structure via the above isomorphisms endows $H^\bullet(X)\oplus H^{\bullet
 -2}(\Delta)(-1)$ with a structure of algebra, which expresses
 cup-product on $\tilde{X}$ as follows.  For $a,b \in
 H^\bullet(X), \; x,y \in H^{\bullet -2}(\Delta)(-1)$ one has:
\begin{enumerate}
\item[] $(0\oplus x)\wedge (0\oplus y)=-j_*(xy)\oplus 2L_\Delta
xy$,

\item[]$(a\oplus0)\wedge(b \oplus0)= ab\oplus 0$,

 \item[] $(a\oplus 0)\wedge (0\oplus y)=0\oplus j^*(a)y$,

  \item[] $(0\oplus x)\wedge (b\oplus 0)=0\oplus xj^*(b).$
\end{enumerate}
\end{enumerate}
The Poincar\'e duality pairing is expressed as follows in terms
  of the above decomposition.  If $x\oplus
  y \in H^i(\tilde{X}), x'\oplus y' \in H^{2n-i}(\tilde{X})$, then
  $$\langle x\oplus y, x'\oplus y'\rangle_{\tilde{X}}=\langle x, x'\rangle_X -\langle y, y'\rangle_\Delta. \blacksquare$$
\end{prop}
  Let $\iota:Y \hookrightarrow X$ denote the canonical inclusion of a
  smooth hyperplane section $Y$ in $X$.  If $Y=X_t$ is a smooth fibre $\rho^{-1}(t)$ of $\rho$, let $k:Y\hookrightarrow
  \tilde{X}$ denote the canonical inclusion.  The following result
  expresses the cohomology of $k_*$ and $k^*$ in terms of
  Proposition \ref{katzuno}.
\begin{prop}\label{katzdos} (Katz ~\cite{SGA7} Exp. XVIII 5.1.1) Notations and assumptions
as above; the restriction homomorphism is expressed by
$$k^*=\iota^*+h_*:H^\bullet(X)\oplus H^{\bullet
-2}(\Delta)(-1)\rightarrow H^\bullet(Y)$$ and the Gysin
homomorphism has the expression
$$k_*=\iota_*\oplus -h^*:H^{\bullet -2}(Y)(-1)\rightarrow H^\bullet(X)\oplus H^{\bullet
-2}(\Delta)(-1).$$

 \end{prop}$\blacksquare$

 Since we will deal with
 the Lefschetz
  theory of both $X$ and $\tilde{X}$, the following discussion will help prove our Main Theorem.

    \begin{caja}{Choice of ${\cal L}_{\tilde{X}}$:} We know by Hartshorne ~\cite{HAG}II.7.10, II.7.11 that the
    line bundle ${\cal L}_N=f^*{\cal L}_X^{\otimes N}\otimes {\cal O}_{\tilde{X}}(-\tilde{\Delta})$ is very ample
    on $\tilde{X}$ for $N\geq N_0.$  Consider $N=m+1$ such that
    $m\geq N_0,$ and choose ${\cal L}_{\tilde{X}}:={\cal
    L}_{m+1}.$
    \end{caja}
\begin{prop}\label{lefschetztilde} Consider the polarisation on $\tilde{X}$ given by the divisor class $\xi_{\tilde{X}}=c_1({\cal
L}_{m+1})=m\cdotp f^*\xi_X+\rho^*([t])$ for $t\in \pp^1$ a regular
value of $\rho$ (not necessary). Let $L_{\tilde{X}}$ be the
Lefschetz operator of this polarisation. One also has
$f^*\xi_X=\xi_X\oplus 0$ and $f^*(\xi_X)\wedge (x\oplus y)=L_X
x\oplus L_\Delta y.$
  In terms of the decomposition of Proposition \ref{katzuno}, $L_{\tilde{X}}^r$ is
   expressed as follows:  $$L_{\tilde{X}}^r(x\oplus 0)=m^{r-1}(m+r)L^r x\oplus -r\cdotp m^{r-1}L_\Delta^{r-1} j^*x,$$
and
   $$L_{\tilde{X}}(0\oplus y)=r\cdotp m^{r-1}L_X^{r-1}j_*y\oplus m^{r-1}(m-r) L_\Delta^r y.$$
\end{prop}
{\bf Proof:} Using Propositions \ref{katzuno}, \ref{katzdos}, the
obtain the following:
 $$\xi_{\tilde{X}}=f^*\xi_X=[Y]\oplus 0, \quad [\tilde{\Delta}]=0\oplus 1_\Delta, \quad [\rho^*(t)]= [Y]\oplus
 -1_\Delta,$$ and $c_1({\cal L}_0)=(m+1)\cdotp [Y]\oplus -1_\Delta.$
 Using the easy fact
 $\rho^*(t)^2=0$,

we obtain
$$\xi_{\tilde{X}}^r=m^r f^*\xi^r_X+r\cdotp
m^{r-1}f^*(\xi^{r-1}_X)\cdotp \rho^*(t)=(m+r)m^{r-1}\xi_X^r\oplus
 -r\cdotp m^{r-1}\xi_\Delta^r.$$  The Proposition now follows from
 Proposition \ref{katzuno}\textit{(ii)}. $\blacksquare$

\section{The Leray filtration of a Lefschetz pencil}
\label{lapiceslefschetz}
    Assume $k=\overline{k}$ as in the previous section.  Choose a Lefschetz pencil on $X$, denoted by $\rho:\tilde{X}\rightarrow\pp^1.$

  Condition {\bf (A)} of Katz ~\cite{SGA7} Exp. XVIII, 5.3 will be important for our
  purposes; we include it below.
\begin{caja}{Condition (A):} Let $\nu:{\cal U}\subset \pp^1$ be contained within the smooth locus of $\rho$.  The adjunction morphisms
 $$R^i\rho_*\Q_\ell \rightarrow \nu_*\nu^*R^i\rho_*\Q_\ell$$
 are isomorphisms for all $0\leq i\leq 2n-2$ (independent of ${\cal U}$).
\end{caja}
An immediate application of the weak Lefschetz theorem yields the first assertion of the following Lemma.
\begin{lemma} \label{lemanu} (~\cite{SGA7} Exp. XVIII Lemma 5.4, Th. 6.3, Cor. 6.4; ~\cite{DeWeII}) If the Lefschetz pencil $\rho$ satisfies condition {\bf (A)}, then the sheaves $R^i\rho_*\Q_\ell$ are constant for $i\neq n-1$.  A suitable multiple of a given polarisation  contains one such Lefschetz pencil.\end{lemma}

\begin{theorem}(\cite{SGA7} 5.6, 5.6.8; ~\cite{DeWeI} Sec. 2; ~\cite{DeWeII})\label{leray} For a Lefschetz pencil $\rho:\tilde{X}\rightarrow\pp^1$, the Leray spectral sequence
$$E_2^{i,j}=H^i(\pp^1,R^j\rho_*\Q_\ell)\Rightarrow H^{i+j}(\tilde{X})$$
degenerates at $E_2$.  For $k:Y=X_t \hookrightarrow \tilde{X}$ the inclusion map of a smooth fibre, the Leray filtration of $\rho$ can be interpreted as follows:
\begin{enumerate}
\item $F^1_\rho H^*(\tilde{X})=\mbox{Ker }k^*$, and $Gr^0_{F_\rho}H^i(\tilde{X})=H^0(\pp^1, R^i\rho_*\Q_\ell)$;
\item $F^2_\rho H^*(\tilde{X})=\mbox{Im }k_*$ image of the Gysin map $k_*$; one has an isomorphism
$F_\rho^2H^i(\tilde{X})=H^2(\pp^1, R^{i-2}\rho_*\Q_\ell)\cong k_*H^{i-2}(Y)$.
\end{enumerate}
The piece $F^2_\rho H^*(\tilde{X})$ coincides with the image of the Gysin homomorphism
$$k_*:H^{*-2}(Y)\rightarrow H^*(\tilde{X}),$$ where $Y$ is a smooth hyperplane section (as above); the piece $F^1_\rho$ coincides with the kernel of $k^*:H^*(\tilde{X})\rightarrow H^*(Y).$
Furthermore, one has $F^2_\rho H^i(\tilde{X})^\perp=F^1_{\rho}
H^{2n-i}(\tilde{X}).$
\end{theorem}
The whole Theorem holds in general.  One can prove the last
assertion under condition {\bf (A)} via the inclusion $F^1_{\rho}
H^{2n-i}(\tilde{X}) \subset F^2_\rho H^i(\tilde{X})^\perp$ and a
dimension count. For the general case see ~\cite{DeWeI} Sec. 2,
~\cite{SGA7}, ~\cite{SGA4}, ~\cite{SGA5} I.5; see also Looijenga
~\cite{Looijenga} Sec. 5. $\blacksquare$

 For the rest of this paper, we fix a Lefschetz pencil on $X$ satisfying condition {\bf (A)}
 above.  We will denote by $L$ the operator in $H^*(\tilde{X})$
 given by $L \bullet=f^*\xi_X\wedge \bullet$ (see Proposition \ref{lefschetztilde}), which has the following expression in terms of \ref{katzuno}:
 \be \label{formulal} L(x\oplus y)=L_Xx\oplus L_\Delta y.\ee
 \begin{remark}  Condition {\bf (A)} induces a Lefschetz theory on the sheaves
  $R^i\rho_*\Q_\ell.$  One has Lefschetz isomorphisms
  $$L^{n-1-i}:R^i\rho_*\Q_\ell \simeq \nu_*\nu^*R^i\rho_*\Q_\ell\longrightarrow\nu_*\nu^*R^{2n-2-i}\rho_*\Q_\ell\simeq R^{2n-2-i}\rho_*\Q_\ell,$$

  where $\nu:{\cal U}\subset \pp^1$ is such that $\rho$ is smooth
  on ${\cal U}.$  We denote by ${\cal P}^i_\rho=\mbox{ker }L^{n-i}$ the primitive cohomology sheaves, and occasionally denote $R^i\rho_*\Q_\ell$ by
  ${\cal R}^i.$
\end{remark}
\begin{corollary}  \label{dualidades}
 \begin{enumerate}
    \item The following isomorphisms hold: $$L^{n-1-i}:R^i\rho_*\Q_\ell \rightarrow R^{2n-2-i}\rho_*\Q_\ell.$$
    \item Let ${\cal P}^i_\rho=\mbox{Ker }L^{n-i}\subset R^i\rho_*\Q_\ell.$  Then
    ${\cal P}^i_\rho$ is constant of fibre $P^i(X)$ if $i\leq n-2$ and
    ${\cal P}^{n-1}_\rho={\cal E}^{n-1}\oplus
    P^{n-1}(X)_{\pp^1},$ where ${\cal E}^{n-1}=\nu_*\nu^*
    {\cal E}^{n-1}$ for all $\nu:{\cal U}\subset \pp^1$ within the
    smooth locus of $\rho$; moreover
    ${\cal E}^{n-1}_t=V(X_t)$ for all geometric
    points
    $\overline{t}\rightarrow t \in {\cal U}.$

    \item Let $0\leq\epsilon\leq 2.$  The pairings
    $$R^i\rho_* \Q_\ell\times R^{2n-2-i}\rho_* \Q_\ell \rightarrow R^{2n-2} \rho_*\Q_\ell\simeq \Q_{\ell}$$
    and
    \be \label{primaris}
     L^{n-1-i}\bullet \cup \bullet:{\cal P}^i_\rho \times {\cal P}^i_\rho \rightarrow \Q_\ell
    \ee
     induce perfect pairings

$$H^\epsilon(R^{i}\rho_*\Q_\ell)\otimes H^{2-\epsilon}(R^{2n-2-i}\rho_*\Q_\ell)\rightarrow H^2(\pp^1,\Q_\ell)$$
and
 \be \label{inducidogr}
  H^\epsilon({\cal P}^i_\rho)\otimes H^{2-\epsilon}({\cal P}^i_\rho)\rightarrow \Q_\ell
\ee

    which agree with the ones resulting from Theorem \ref{leray}; for instance, the pairing given by
    $a\otimes b \mapsto \langle L^{n-1-i} a, b\rangle_{\tilde{X}}$ in $\mbox{Gr}_{F_\rho}^{\bullet} H^*(\tilde{X})$ equals
    the one in (\ref{primaris}).

    One has $\mbox{dim }H^0({\cal P}^i_\rho)=\mbox{dim }H^2({\cal
    P}^i_\rho)$ for $0\leq i\leq n-1$.
\item The Lefschetz isomorphisms on sheaves translate also into
 their cohomology groups; in particular, $$H^\epsilon({\cal
P}^i_\rho)=\mbox{ker }(L^{n-i}:H^\epsilon({\cal R}^i)\rightarrow
 H^\epsilon({\cal R}^{2n-i})).$$
 \item $\mbox{dim
}H^0({\cal R}^i)=\mbox{dim }H^2({\cal R}^i)=b_i(X)$ for all $i$.
As a result, $\mbox{dim }H^0({\cal P}_\rho^i)=\mbox{dim }H^2({\cal
P}_\rho^i)=\mbox{dim }P^i(X).$
\end{enumerate}

\end{corollary}

{\bf Proof:} The result follows from Lemma \ref{lemanu}, Theorem
\ref{leray}, Deligne ~\cite{DeWeI} 2.8 and 2.12, and Katz
~\cite{SGA7} XVIII Lemmas 5.4, 5.5, 5.6.9 and proof of Th. 5.6.8.

Let us check the last assertion for $i=n-1:$  the morphism
$k_*=\iota_*\oplus -h^*:H^{n-1}(Y)$ has kernel $V(Y)=\mbox{ker
}\iota_*.$  Therefore $\mbox{dim }H^2({\cal
 R}^{n-1})=b_{n-1}(X)=b_{n-1}(Y)-\mbox{dim }V(Y).$ The equality
${\cal R}^i={\cal P}^i_\rho \oplus L{\cal R}^{i-2}$ yields
$\mbox{dim }H^0({\cal P}^i_\rho)=\mbox{dim }P^i(X).$
(Alternatively, use ~\cite{SGA7} XVIII Th. 5.6.)
  $\blacksquare$

  We now view the computations of Proposition
 \ref{lefschetztilde} in a different fashion.
\begin{lemma}\label{calculostildes2} Let $x\oplus y \in H^i(\tilde{X})$, and let $r\in \N.$
Then $\xi_{\tilde{X}}-m\cdotp
 f^*(\xi_X)=\rho^*([t])=k_*(1_{H^*(Y)}) \in F^2_\rho.$  Thus the
expression
 \be\label{diferenciaeles} (L_{\tilde{X}}^r-m^rL^r)(x\oplus y)=L_{\tilde{X}}^r(x\oplus y)-m^r (L^rx\oplus L^r_\Delta y)=r\cdotp m^{r-1} k_*L^{r-1}_Y(\iota^*x+h_*y) \ee
belongs to $F^2_\rho.$
\end{lemma}
{\bf Proof:} By (\ref{formulal}) we have $L^s(x\oplus
y)=f^*\xi_X^s \wedge (x\oplus y)=L^sx\oplus L_\Delta^sy.$ On the
other hand, if $Y=X_t$ is a smooth geometric fibre, then
$\rho^*([pt.])=k_*(1_Y)=\xi_X\oplus -1_\Delta \in H^2(\tilde{X}).$
We have $\rho^*(t)\wedge (x\oplus 0)=Lx\oplus -j^*x=k_*(\iota^*x)$
and $\rho^*(t)\wedge (0\oplus y)=j_*y\oplus -L_\Delta
y=k_*(h_*y)$, whence $\rho^*(t)\wedge (x\oplus
y)=k_*(\iota^*x+h_*y)=Lx+j_*y\oplus -(j^*x+L_\Delta y).$ Finally
\bean r\cdotp m^{r-1} L^rx+L^{r-1}j_*y\oplus -L^{r-1}_\Delta(
j^*x+L_\Delta y)=L_{\tilde{X}}^r(x\oplus y)-m^r (L^{r}x\oplus
L^{r}_\Delta y)=\\
=r\cdotp m^{r-1}k_*L^{r-1}_Y(\iota^*x+h_*y)\eean as desired.
$\blacksquare$

\begin{corollary}\label{lefschetztildedec} Notations and assumptions
being as above,
$$L^{n-i}_{\tilde{X}}(P^i(X)\oplus 0)=L^{n-i}P^i(X)\oplus
0=k_*L^{n-i-1}_Y \iota^*P^i(X)\subset F^2_\rho$$ and
$P^i(\tilde{X})\supset P^i(X)\oplus 0.$  One has
 $L^r_{\tilde{X}}k_*(y)=m^rk_*(L^r_Yy)=m^rL^rk_*y$ for all $r\geq 0.$

\end{corollary}
{\bf Proof:}  By ~\cite{Kleiman} 1.4.7, $h^*L^{n-1-i}_Y P^i(Y)=0$
and
$$k_*:L^{n-i-1}_Y \iota^*P^i(X) \rightarrow L^{n-i}P^i(X)\oplus 0$$ is an
isomorphism.  Let us prove the inclusion $P^i(X)\oplus 0\subset
P^i(\tilde{X}).$  By formula (\ref{diferenciaeles}), it suffices
to check that
$$(n-i)m^{n-i}k_*L_Y^{n-i-1}P^i(Y) \subset L^{n-i}P^i(X)\oplus
0,$$ but this inclusion is clear.  We have seen that the image of
$P^i(X)\oplus 0$ via the Lefschetz isomorphism is precisely
$L^{n-i}P^i(X)\oplus 0,$ thus establishing the result.
$\blacksquare$
\begin{remark} By Lemma \ref{calculostildes2}, the operator $L_{\tilde{X}}-m\cdotp L$ vanishes on
$Gr^*_FH^*(\tilde{X}).$ The same thing happens on the sheaves
${\cal R}^i$.
\end{remark}
\begin{corollary} \label{primisyceros} The map $L^{n-i}$ and the Lefschetz isomorphism
$L^{n-i}_{\tilde{X}}$ yield isomorphisms $$(P^i(X)\oplus 0)\oplus
 k_*H^{i-2}(Y) \overset{\sim}{\rightarrow} k_*H^{2n-2-i}(Y).$$ The
subspace $L_{\tilde{X}}^{j}P^i(X)\oplus 0$ is linearly disjoint
 with $F^1_\rho$ for $j<n-i$, and $L^{n-i}P^i(X)\oplus 0\subset
 F^2_\rho.$
\end{corollary}
 The first assertion follows from Corollary \ref{lefschetztildedec} and Corollary \ref{dualidades}\textit{(5)}.  The second
  assertion follows from the first. $\blacksquare$

\begin{corollary}\label{inicioprim} The natural map $P^i(X)\oplus 0 \rightarrow H^0({\cal R}^i)$ of Theorem \ref{leray} induces an isomorphism
$$P^i(X)\oplus 0 \cong H^0({\cal P}^i_\rho)$$ for $0\leq i\leq
n-1.$  The map $\rho^*(t)\wedge \bullet$ yields an isomorphism
between $H^0({\cal P}^i_\rho)$ and $H^2({\cal P}^i_\rho)$.  As a
result, $H^2({\cal P}^i_\rho)=LP^i(X)\oplus P^{i-2}(\Delta)\cap
F^2_\rho H^{i+2}(\tilde{X})=k_*P^i(Y)$ for $i\leq n-1.$
\end{corollary}
{\bf Proof:}
\begin{enumerate}
\item The dimensions are equal, and $L^{n-i}(P^i(X)\oplus
0)\subset F^2_\rho,$ hence the map $$P^i(X)\oplus 0 \rightarrow
H^0({\cal R}^i)$$ induces an isomorphism onto $H^0({\cal
 P}^i_\rho).$
\item The class $[\rho^*(t)] \in \rho^*H^2(\pp^1)$, hence
$\rho^*(t)\wedge \bullet$ induces a map $H^0({\cal
P}^i_\rho)\rightarrow H^2({\cal P}^i_\rho),$ which reads as
follows: $$\rho^*(t)\wedge (x\oplus 0)=k_*k^*(x\oplus
0)=k_*\iota^*x.$$  Therefore its image is $k_*\iota^*P^i(X),$
 whose dimension agrees with $\mbox{dim }H^2({\cal P}^i_\rho).$
The assertion is thus proven.
 \end{enumerate}
 $\blacksquare$

\subsection{Absolute and relative correspondences}\label{relcorresps}
   Let $p:M\rightarrow B$ be a smooth projective morphism onto a
smooth algebraic variety $B$; denote the dimension of $M$ by $n$.
In this section we will establish the usual properties of the
composition law of correspondences on $M\times_B M$; if $u$ is a
codimension-$(r-\mbox{dim }B)$ cycle on $M\times_B M$, the degree
of $u$ as a relative correspondence is defined to be $r$, i.e. the
same as that of $u$ as a correspondence of $M$; for instance, this
definition makes the cycle $\Delta_M$ into a relative
correspondence of degree $0$.  See Fulton ~\cite{F} Ch. 10, 16 for
an introduction.  The heart of this section is Lemma
\ref{printer}, where every identity holds modulo rational
equivalence.  If $u,v \in CH^{n-1+*}(M\times_BM)$, then the
composition of $u, v$ relative to $B$ is defined to be
$$v\circ_Bu:=p_{13*}^B(p_{12}^{B*}(u)\bullet p_{23}^{B*}(v)),$$
where $p^B_{ij}:M\times_B M\times_B M \rightarrow M\times_B M$ are
the canonical projections.  $\circ_B$ endows $CH^{n-1+*}(M\times_B
M)$ with a ring structure, and the usual properties hold.  The
upshot is Proposition \ref{fil}.
\begin{lemma}\label{printer} Notations and assumptions as above: suppose that $$u, v\in CH^{*}(M\times_B M)$$ are relative
correspondences of degrees $r,s$ respectively.  If $t\in B$ is a
closed point, let $\lambda_t:M_t\times M_t \hookrightarrow
M\times_BM$ be the canonical inclusion, and let
$u_t=\lambda_t^*(u).$  Then:\begin{enumerate} \item The
correspondence $v\circ_B u$ is of degree $r+s$; \item for any
$t\in B$ one has: $(v\circ_B u)_t=v_t\circ u_t.$ \item one can
compare the composition laws $\circ$ and $\circ_B$ as follows. Let
$$j:M\times_BM\hookrightarrow M\times M$$ denote the natural
inclusion; then $j_*(v\circ_B u)=j_*v\circ j_*u.$
\end{enumerate}
\end{lemma}
{\bf Proof:} The first and second assertions are clear.  For the third assertion, we need some notation:
 Let $p^B_{ab}$ denote the $(a,b)$-projections $M\times_B M\times_B M\rightarrow M\times_B M$ and $p_{ab}$ denote the corresponding projections $M^3\rightarrow M^2$.  Let
$$inc:M\times_B M\times_B M \hookrightarrow M\times M\times M$$
be the natural inclusion.  For each pair $a\neq b$ in $\{1,2,3\}$ we have a fibre product
\be
\begin{CD}\label{prodfibr}
 N_{ab} @>j'_{ab}>>  M^3  \\
            @Vp_{ab}^{\prime}VV  @VVp_{ab}V \\
 M\times_B M @>j>> M^2.
\end{CD}
\ee

 Denote by $k_{ab}:M\times_B M\times_B M \hookrightarrow N_{ab}$ the natural inclusion.  Then $p'_{ab}k_{ab}=p_{ab}^B$.

 We wish to prove the identity
$$j_*p_{13*}^B(U\times_B M\bullet M\times_B V)=p_{13*}(U\times M\bullet M\times V);$$
by the above one has $j_*p_{13*}^B=p_{13*}inc_*$.  Thus it suffices to prove the following:
$$inc_*(U\times_B M\bullet M\times_B V)=U\times M\bullet M\times V.$$
In other words, one must check, for $u,v \in CH^*(M\times_B M)$:
\be \label{ident} p_{12}^*(j_*u)\bullet
p_{23}^*(j_*v)=inc_*(p_{12}^{B*}(u)\bullet p_{23}^{B*}(v)). \ee

Consider now the following cartesian diagram of embeddings: \be
\label{fibra}
\begin{CD}
M\times_B M\times_B M  @>inc>> M^3 \\
 @V(k_{12}, k_{23})VV  @VV\Delta_{M^3}V \\
N_{12}\times N_{23}  @>>j'_{12}\times j'_{23}> M^3\times M^3.
\end{CD}
\ee Note that \be\label{diagoskas} (k_{12}, k_{23})=(k_{12}\times
k_{23}) \Delta_{M\times_B M\times_B M}. \ee  Formula (\ref{ident})
is equivalent to the following: \be \label{idfinal} \Delta_{M^3}^*
(p_{12} \times p_{23})^* (j_*\times j_*)= inc_*\Delta_{M\times_B
M\times_B M}^*(p_{12}^B\times p_{23}^B)^*. \ee

  Let us develop the l.h.s. Using Fulton ~\cite{F} Prop. 1.7, one has $p_{ab}^*j_*=j'_{ab*} p_{ab}'^*$ from the fibre product (\ref{prodfibr}), hence
$$\mbox{l.h.s.}=\Delta_{M^3}^*(j'_{13}\times j'_{23})_*(p'_{12}\times p'_{23})^*.$$  Since
 all varieties involved are smooth and quasiprojective, (\ref{fibra})
yields \be\label{moviendo} \Delta_{M^3}^*(j'_{13}\times
j'_{23})_*=inc_*(k_{12},k_{23})^*:\ee
  Indeed, by the moving lemma ~\cite{Roberts} one can check the
  above for an algebraic cycle $\zeta$ in $N_{12}\times N_{23}$ that intersects
  properly with the image of $(k_{12}\times k_{23}).$  The identity (\ref{moviendo}) is easily
  established once for such $\zeta$, and so for each Chow class in $N_{12}\times N_{23}.$
 Using (\ref{diagoskas}) and (\ref{moviendo}) yields
\be\label{conyoliso} \mbox{l.h.s.}=inc_* \Delta_{M\times_B
M\times_B M}^*(k_{12}\times k_{23})^*(p'_{12}\times p'_{23})^*,
\ee
 which in turn equals
$$inc_*\Delta_{M\times_B M\times_B M}^*(p^B_{12}\times p^B_{23})^*=\mbox{r.h.s};$$
the Lemma is thus established. $\blacksquare$

\begin{remark} Lemma \ref{printer} generalises acccordingly when source and target of the relative correspondences are different; so do Lemma \ref{preserv}, Proposition \ref{fil}.\end{remark}
\begin{lemma}\label{preserv} With the notations and assumptions of Lemma \ref{printer}, we denote by \newline
$\iota=\iota_t:M_t \hookrightarrow M$ the canonical inclusion.  The following identity holds:
$$j_*u \circ \iota_*=\iota_*\circ u_t.$$

\end{lemma}
{\bf Proof of Lemma \ref{preserv}:}  Consider the following
cartesian diagram of embeddings:
 \bean
 \begin{CD}
M_t\times M_t  @>1\times\iota>> M_t\times M \\
 @V\lambda_tVV  @VV\iota\times 1V \\
M\times_BM  @>j>> M\times M.
\end{CD}
\eean
  Analogously as shown in Lemma \ref{printer} with (\ref{fibra}), the following formula
  holds:
\be \label{trucofibra}(\iota\times 1)^*j_*=(1\times
\iota)_*\lambda_t^*.\ee
  Using Fulton ~\cite{F} 16.1.1.(c) (see also Scholl ~\cite{Scholl}
  1.10) we derive $$j_*(u)\circ \iota_*=\iota_*\circ u_t.$$  The
  proof is now complete.
$\blacksquare$

\begin{prop} \label{fil} The correspondences supported on $D=\tilde{X}\times_{\pp^1}\tilde{X}$ preserve the Leray filtration.  More precisely, if $u \in
CH_*(\tilde{X}\times_{\pp^1}\tilde{X})$ and
$j':\tilde{X}\times_{\pp^1}\tilde{X}\hookrightarrow
\tilde{X}\times\tilde{X}$, then $[j'_*(u)]F_\rho^i\subset
F_\rho^i$ for $i=0,1,2$. In addition, if $u$ is
 supported on a finite set of fibres of the structure morphism $\tilde{X}\times_{\pp^1}\tilde{X}\rightarrow\pp^1$, then
 $[j'_*u]|F^2_\rho=0, [j'_*u]|H^*(\tilde{X})\subset F^1_\rho.$\end{prop}
{\bf Proof of Proposition \ref{fil}:}

   Choose a smooth fibre $Y=\tilde{X}_t$ of $\rho$.  We will establish an identity identical to that of Lemma
   \ref{preserv}, proved with the due care since $D$ is not smooth in general.  Let $B\subset \pp^1$ denote
the smooth locus of
   $\rho$ and let $D_B:=\rho^{-1}(B)\times_B
\rho^{-1}(B)$ be the smooth locus of $D/\pp^1.$  Now $\lambda_t$
factors as
$$X_t\times X_t \overset{\lambda'_t}{\hookrightarrow} D_B
\overset{\nu'}{\hookrightarrow} D,$$ where $\nu'$ is an open
immersion and $D_B$ is smooth.  We may now define
$\lambda_t^*=\lambda_t^{\prime *}\nu'^*.$  The identity
$(\iota\times 1)^*j'_*z=(1\times \iota)_*\lambda_t^*z$ holds for
any algebraic cycle on $D$: if $z$ is supported on $X_t \times
X_t$, then both sides are $0$ by Fulton ~\cite{F} 10.1 (use
\emph{op.cit.} Cor. 6.3).  If $z$ has no component contained in
$X_t \times X_t$, then $z$ intersects properly with $X_t \times
X_t$ and $j_*z$ intersects properly with $X_t\times \tilde{X}$;
again, it is easy to check that both sides agree (as algebraic
cycles, no equivalence relation established). Just as in Lemma
\ref{preserv}, then for every correspondence $j'_*u$ on
$\tilde{X}$ supported in $D$ we have
 \be \label{filtrosrho} j'_*u\circ \iota_*= \iota_*\circ u_t.\ee
 Now we check that any correspondence supported on
   $D=\tilde{X}\times_{\pp^1}\tilde{X}$ sends $F^2$ into $F^2$: indeed, (\ref{filtrosrho}) directly implies $j'_*(u)F^2_\rho\subset F^2_\rho.$
     Likewise, $j'_*(u)$  sends $F^1_\rho$ into $F^1_\rho$: suppose $x\in H^*(\tilde{X})$ is such that $\iota^*x=0.$  We have $\iota^*\circ j'_*(u)= u_t
     \iota^*$,
      hence $j'_*(u)F^1 \subset F^1.$  In the case when $u$
   is supported on a finite union of subschemes $X_{s_i}\times X_{s_i}$, taking $t\neq s_1, \cdots, s_r$ in the above argument
   yields $u_t=0$, thus proving the second assertion.  The proof is now complete.  $\blacksquare$

\subsection{More on supports}

  We now return to the setting of Theorem \ref{leray} and assume the notations of Theorem \ref{leray} and Proposition \ref{fil}.  

Let $u, v$ be
  algebraic cycles supported on $D$.  If $p^*_{12}(u), p^*_{23}(v)$
  intersect properly, then the correspondence $v\circ u$ is easily
  seen to be supported on $D$; a similar argument works on Chow classes if an embedded desingularisation to
  $\tilde{X}\times_{\pp^1}\tilde{X}\subset \tilde{X}^2$ exists --
  in the proof of Lemma \ref{printer}, formula
  (\ref{conyoliso}) required smoothness of the integral scheme $M\times_BM$;
   In a similar vein to \ref{relcorresps}, one derives the
following statement.
\begin{prop}\label{soportes} Notations and assumptions being as above,  let $u=j'_*u_0,v=j'_*v_0$ be
two correspondences supported on $D$.  Let $j''_*:M\times_BM
\subset M\times M$ denote the canonical inclusion.

  The following
statements hold (modulo rational equivalence).
\begin{enumerate}
    \item Let $r,r_0$ denote the inclusions $r:M\times M\subset
    \tilde{X}\times \tilde{X}$, $r_0:M\times_BM\subset D$. Then $r^*v\circ r^*u=j''_*(r_0^*v_0\circ_B
    r_0^*u_0).$
    \item Let $\phi:Z\twoheadrightarrow D\times \tilde{X}$ be a De
    Jong alteration ~\cite{dJ}.  Let $z$ be such that
    $\phi_*(z)=p_{12}^*(u).$ Then $p_{12}^*(u)p_{13}^*(v)$ is supported on
    $\tilde{X}\times_{\pp^1}\tilde{X}\times_{\pp^1}\tilde{X}$;
    as a result, $v\circ u$ is supported on $D$.
\end{enumerate}
\end{prop}
{\bf Proof:} The first statement follows from Lemma \ref{printer}.
Let us prove the second statement: for $u, z, \phi$ as above, one
has $\phi_*(z\bullet \phi^*(p_{23}^*(v))=\phi_*(z'\bullet
\phi^*(p_{23}^*(v))$, where $z'\sim_{\mathrm{rat}} z$ is such that
$z'$, $\phi^*(p_{23}^*(v))$ intersect properly.  It is now
apparent that the Chow class
$$\phi_*(z'\bullet \phi^*(p_{23}^*(v))=\phi_*(z')\bullet p_{23}^*(v)= p_{12}^*(u)\bullet p_{23}^*(v)$$
is
 supported on
$\tilde{X}\times_{\pp^1}\tilde{X}\times_{\pp^1}\tilde{X}$.
Applying $p_{13*}$ yields the Chow class $v\circ u$, which is
thererefore supported on $D$, thus completing the proof.
$\blacksquare$

\subsection{Action on the Leray spectral sequence}
\label{accionss}
  (Again we assume $k$ algebraically closed.)  Let $u$ be a correspondence of degree $r$ supported on
  $\tilde{X}\times_{\pp^1}\tilde{X}$.  Then $u$ induces a
  correspondence $u_t$ of degree $r$ on $X_t$ for each $t\in
  B(\overline{k})$, where $\nu:B\rightarrow \pp^1$ is the smooth
  locus of $\rho$ as above.  $u$ thus defines a homomorphism of $\ell$-adic sheaves for $0\leq j\leq n-1$:
\be\label{mfladicos}
  u:\nu_*\nu^*R^j\rho_*\Q_\ell=R^j\rho_*\Q_\ell\rightarrow
  \nu_*\nu^*R^{j+2r}\rho_*\Q_\ell=R^{j+2r}\rho_*\Q_\ell\ee
  (using {\bf (A)}), which in turn yields $\Q_{\ell}$-linear maps
  \be\label{indcorrespfibras} H^i(R^j\rho_*\Q_\ell)\rightarrow H^i(R^{j+2r}\rho_*\Q_\ell).\ee
These maps clearly agree with those induced on
  $Gr^*_{F_\rho}$ by $j_*u$ in Proposition \ref{fil}, and so do the respective composition laws.
\begin{obs}\label{remarcagenerico} Morphisms induced in (\ref{mfladicos}) and (\ref{indcorrespfibras}) depend only on the class of $u$ in
$H^*({\cal Y}\times {\cal Y})$.  Indeed, denote the generic point
of $\pp^1$ by $\eta$, and the image of $u$ in $A^{n-1+*}({\cal
Y}\times {\cal Y})$ by $u_\eta$ or $[u]_{\cal Y}$; suppose that
$u'_\eta-u''_\eta|H^j({\cal Y})=0$ for all $j$. Then for a
sufficiently small neighbourhood $\nu_1:{\cal U}_1\subset \pp^1$
of $\eta$ one has
$0=u'-u''|\nu_{1*}\nu_1^*R^j\rho_*\Q_\ell=R^j\rho_*\Q_\ell$ for
all $j$, hence $u'-u''$ induces $0$ on
$Gr^*_{F_\rho}H^*(\tilde{X}).$
 Here we used ~\cite{FKiehl} I.12.10, I.12.13 (see also ~\cite{SGA5}) and the base change theorems in etale cohomology
  ~\cite{SGA4}; ~\cite{FKiehl} I.6,I.7.\end{obs}
\begin{caja}{Definition.} Let ${\cal A}\subset A^{n+*}(\tilde{X}\times\tilde{X})$ denote the subring (see Proposition \ref{soportes})
 of homological correspondences
supported on $\tilde{X}\times_{\pp^1}\tilde{X}.$ Let ${\cal J}$ be
the ideal of \newline $A^{n+*}(\tilde{X}\times\tilde{X})$
consisting of the elements $u$ such that $uF^i_\rho \subset
F^{i+1}_\rho$ for $i=0,1,2$ (i.e. those inducing $0$ on
$\mbox{Gr}^*_{F_\rho}H^*(\tilde{X})$).  Let ${\cal I}$ be the
ideal (see Lemma \ref{cardosprevio}) of ${\cal A}$ consisting of
the $u\in {\cal A}$ such that $u=[j'_*v]$ with $v$ an algebraic
cycle on $\tilde{X}\times_{\pp^1}\tilde{X}$ (with
$\Q$-coefficients) inducing an homologically trivial class on
$H^*({\cal Y}\times {\cal Y}).$ One has ${\cal I}\subset {\cal J}$
by (\ref{mfladicos}).  We denote by ${\cal K}\subset {\cal A}$ the
subspace of all classes $w=[j'_*w_0] \in {\cal A}$
    satisfying $[w_0]_{\cal Y}H^*{\cal
    Y} \subset V({\cal Y}).$
\end{caja}
  The following Proposition sharpens Remark \ref{remarcagenerico} above.
\begin{prop}\label{generico} Then the ideals ${\cal I},{\cal J},{\cal K}$ of ${\cal A}$
satisfy ${\cal I}\subset {\cal J}, {\cal I}^{\circ 3}={\cal
J}^{\circ 3}=0.$  Let $w_0$ be an algebraic cycle supported on
$D$, representing the correspondence $w\in
A^{n+r}(\tilde{X}\times\tilde{X})$
 of degree $r$.
Suppose that $r\neq 0$; then $w_0 \in {\cal I}$

  if and only if $w\in {\cal K}.$
In general, the following statements hold.
 \begin{enumerate}
    \item[(i)] ${\cal J}\subset {\cal K}.$  ${\cal K}$ is an ideal of ${\cal A}.$ 
    \item[(ii)] Suppose that $H^1(R^{n-1}\rho_*\Q_\ell)=0.$  Then ${\cal J}={\cal K}.$
    \item[(iii)] Assume that $n$ is even or $\mbox{char }k\neq 2$.  If $H^1(R^{n-1}\rho_*\Q_\ell)\neq
     0$, then ${\cal I}={\cal J} \varsubsetneq {\cal K}.$
    \item[(iv)] If $n$ is even or $\mbox{char }k\neq 2$, then there exists $d_0\in \N$ such that, for every $d\geq d_0$, every Lefschetz fibration of degree-$d$
    hypersurfaces satisfies \textit{(iii)}.
    \item[(v)]  If $n$ is even or $\mbox{char }k\neq 2$, and $H^1(R^{n-1}\rho_*\Q_\ell)\neq 0$,
    then any Chow class $u$ supported on $D$ such that $[j'_*u]=0$ satisfies
 $[u]_{\cal Y}=0.$

\end{enumerate}
\end{prop}

{\bf Proof:}  The nilpotence assertion for ${\cal I}, {\cal J}$ is
clear.  Now, let $w_0$ be as above with $r\neq 0;$ we argue as in
Remark \ref{remarcagenerico}.  If ${\cal R}^{i}$ is constant
($i\neq n-1$) then  $[w_0]_{{\cal Y}} H^i({\cal Y})=0$ if and only
if $w_0:{\cal R}^{i}\rightarrow {\cal R}^{i+2r}$ is $0.$ The same
holds whenever ${\cal R}^{i+2r}$ is constant. Thus $w\in {\cal
I}\Leftrightarrow w\in {\cal J}$ whenever $r\neq 0.$

  Assume $r=0$.  Then ${\cal R}^{n-1}=L{\cal R}^{n-3}\oplus {\cal P}^{n-1}_\rho$, where
  ${\cal P}^{n-1}_\rho=P^{n-1}(X)_{\pp^1} \oplus {\cal E}^{n-1}$, and $V({\cal Y})={\cal
  E}^{n-1}_{\overline{\eta}}$ is the geometric
  generic fibre of ${\cal E}^{n-1}$ (here $\overline{\eta}$ is a geometric generic point of
  $\pp^1.$) Hence $w\in {\cal J} \Leftrightarrow w_0 {\cal
  R}^\bullet= w_0 {\cal E}^{n-1}
  \subset {\cal E}^{n-1}.$  The $\Q_\ell$-adic sheaf $w_0{\cal
  E}^{n-1}$ has $[w_0]_{\cal Y}V({\cal Y})$ as its geometric generic
  fibre; all the pieces of $H^*({\cal Y})$ are monodromy invariant, except $V({\cal Y})$ which has no invariants (by {\bf (A)}).  This settles \emph{(i)},\emph{(ii)}.

    To prove \emph{(iii)}, recall that if $n$ is even or
    $\mbox{char }k\neq 2$, then all the singularities of $\rho$
    are non-degenerate quadratic singularities of fibres, and the monodromy
    representation of $\pi_1^{alg}(B, \overline{\eta})$ ($B$ being the smooth locus of $\rho$) on $V({\cal Y})$ is
    absolutely irreducible (~\cite{SGA7}, esp. XVIII Cor. 6.7).
    As a result, the $\pi_1(B)$-submodule $[w_0]_{\cal Y}V({\cal Y})$ is either
    $0$ or $V({\cal Y}),$ and so there are two possibilities for
 the inclusion of $\Q_\ell$-sheaves $w_0{\cal
  E}^{n-1}\hookrightarrow {\cal E}^{n-1}$: either the image or the cokernel of
  this inclusion are skyscraper sheaves.  Since $H^1({\cal R}^{n-1})=H^1({\cal E}^{n-1})\neq 0$, we
    have $$wH^1({\cal R}^{n-1})=0\Leftrightarrow [w_0]_{\cal Y}V({\cal
    Y})=0,$$ thus settling \textit{(ii)},\emph{(iii).}

   Let us prove \emph{(iv)}; by Lemma \ref{piezasola} (see ~\cite{SGA7} XVIII Th. 5.7), $H^1({\cal R}^{n-1})\simeq P^n(X)\oplus V(\Delta).$
    If $P^n(X)\neq 0$ there is
  nothing to prove; if $P^n(X)=0$ the assertion follows from the
  next elementary lemma.
  \begin{lemma} (compare ~\cite{SGA7} XVIII Lemme 6.4.2) \label{bettid} With the notations and hypotheses of this section (assuming $n\leq 3$),  let $d\in \N$.  Let $Y(d), Y'(d)$ denote degree-$d$ hypersurface
  sections intersecting transversally, and let $\Delta(d)=Y(d)\cap Y'(d).$  Then $b_{n-2}(\Delta(d))$ is a polynomial of degree $n$ in $d$.\end{lemma}
{\bf Proof of Lemma \ref{bettid}:} Let $c(X), c(\Delta)$ be the
total Chern classes of $X, \Delta$ and let $j:\Delta \subset X$
denote the canonical inclusion. 
  Let $\int_X$ denote the trace map on $X$, and $H=c_1({\cal
O}_X(1))$ with the present polarisation.  Then, using
$j_*j^*\ga=d^2H^2\bullet \ga$, we obtain:
$$\chi(\Delta(d))=\int_{\Delta(d)} c(\Delta(d))=\int_{\Delta(d)} j^*\frac{c(X)}{(1+d\cdotp
H)^2}=\int_X \frac{d^2H^2c(X)}{(1+d\cdotp H)^2},$$ which is a
polynomial in $d$ of degree $n$ with lead term $(-1)^{n} \mbox{deg
}X\cdotp d^n.$  Isolating yields
$b_{n-2}(\Delta(d))=(-1)^{n-2}\chi(\Delta(d))+2\sum_{i\geq
1}(-1)^i b_{n-2-i}(X)$, thus completing the proof. $\blacksquare$

  Taking $d\gg 0$, the $d$-uple embedding of $X\subset \pp$ satisfies the hypotheses of \emph{(iii)}, thus establishing
  \emph{(iv)}.

    It remains to prove \textit{(v)}. By Corollary \ref{dualidades}, $H^1({\cal E}^{n-1})=H^1({\cal R}^{n-1})$. 

   If $[j'_*(u)]=0 \in {\cal J}$, then $u$ induces $0$ on $H^\epsilon({\cal
    R}^i)$ for all $\epsilon, i$ and from  \textit{(iii)} we derive $[u]_{\cal Y}=0$.

    Proposition \ref{generico}
  is thus established. $\blacksquare$

\begin{corollary}\label{cocientebien} If $n$ is even or $\mbox{char }k\neq 2$ and $H^1({\cal R}^{n-1})\neq 0$,  then the
restriction map
$$CH_{n-1-*}(D)\twoheadrightarrow A^{n-1+*}({\cal Y}\times {\cal Y})$$ factors through a ring homomorphism
$$\mbox{res}_{\cal Y}:{\cal A}\twoheadrightarrow A^{n-1+*}({\cal Y}\times {\cal Y})$$
whose kernel is $\mbox{res}_{\cal Y}={\cal I}.$
\end{corollary}
{\bf Proof:}  The Corollary follows from Proposition
\ref{generico}\textit{(v)}.  $\blacksquare$

  The next Lemma is only necessary if $n$ is odd and $\mbox{char }k=2.$
\begin{lemma}\label{cardosprevio} Notations and assumptions as above.  The subspace ${\cal I}\subset {\cal A}$ is an ideal.\end{lemma}
{\bf Proof:}  Let $B\subset \pp^1$ be the smooth locus of $\rho$,
and $M=\rho^{-1}(B)$. We define the map
$[j'_*]:CH_{n-1-\bullet}(D)\rightarrow A^{
n+\bullet}(\tilde{X}\times \tilde{X})$ (here $[\bullet]$ is the
cycle map, and $D=\tilde{X}\times_{\pp^1}\tilde{X}$) whose image
is precisely ${\cal A}$.  If $r, r_0$ are as in Proposition
\ref{soportes} and \newline $j'':M\times_B M \hookrightarrow
M\times M$, we have a commutative diagram

\bean
 \begin{CD}
  CH_{n-1-*}(D)@>[j'_*]>> A^{n+*}(\tilde{X}\times \tilde{X}) \\
 @V[r_0^*]VV  @VVr^*V \\
A^{n-1+*}(M\times_B M)@>j''_*>> A^{n+*}(M\times M).
\end{CD}
\eean
  We have $r^*{\cal A}=j''_*A^{n-1+*}(M\times_B M)$, and $r^*, j''_*$
  are
  ring homomorphisms by Proposition \ref{soportes}.  It is clear
  that $r^*{\cal I}$ is an ideal of $r^*{\cal A},$ which coincides
  with the image of $\mbox{Ker}(A^{n-1+*}(M\times_B M) \twoheadrightarrow A^{n-1+*}({\cal Y}\times {\cal
  Y}))$ via $j''_*.$  Now the kernel of ${\cal I}\twoheadrightarrow
  r^*{\cal I}$ consists of the correspondences supported on $\bigcup_{s\in \pp^1 \mathrm{\;singular}} X_s\times
  X_s.$ It is now clear that ${\cal I}\subset {\cal A}$ is an ideal. $\blacksquare$

\begin{obs}\label{remarcapostgenerico}
  By Proposition \ref{generico} above, there is a ring epimorphism
$$\varphi=\varphi_{\cal Y}:A^{n-1+*}({\cal Y}\times{\cal Y}) \twoheadrightarrow {\cal A}/{\cal I},$$
which is an isomorphism if $n$ is even or $\mbox{char }k\neq 2$ by
Corollary \ref{cocientebien}.  It is not clear whether
$\varphi_{\cal Y}$ is an isomorphism if $n$ is odd and $\mbox{char
}k=2.$
\end{obs}

  The next Corollary circumvents the possible non-isomorphy of
  $\varphi$ for the purposes of this paper.  Its proof is
  straightforward.
\begin{corollary}\label{rhocardos} Let ${\mathfrak a}\subset A^{n-1+*}({\cal Y}\times{\cal Y})$ be the ideal of correspondences
 $u$ such that $uH^*({\cal Y})\subset V({\cal Y}).$ One has a ring isomorphism induced by $\varphi$ above:
$$A^{n-1+*}({\cal Y}\times{\cal Y})/{\mathfrak
a} \overset{\sim}{\rightarrow} {\cal A}/{\cal K}.$$
 Consider a graded unital subalgebra ${\cal B}\subset A^{n-1+*}({\cal Y}\times {\cal Y})$
such that ${\cal B}\cap {\mathfrak a}=0.$
 Then $\varphi$ yields an isomorphism ${\cal B}\simeq
\varphi({\cal B})\subset {\cal A}/{\cal I},$ which maps
isomorphically after composing with the quotient map ${\cal
A}/{\cal I}\rightarrow {\cal A}/{\cal K}.$
\end{corollary}
$\blacksquare$

\section{The relative projectors}
\label{proyrelativos}

  We have seen in Lemma \ref{adjuntoh} that, if $C(X)$ holds, then the ring of
  correspondences of $X$, $A^{\mathrm{dim }X+\bullet}(X\times X)$
  decomposes through the adjoint action of $H_X$, $u\mapsto [H,u]$; the
  degree-$0$ correspondences are exactly those commuting with
  $H_X$, or equivalently, with the K\"unneth projectors $\pi^i_X$ for all $i$.  We wish to translate this
  situation into the relative context presented in Section
  \ref{lapiceslefschetz}.  Our first goal is to create natural relative
  analogues $\pi^i_\rho, H_\rho$
  of $\pi^i$ and of $H$, supported on $\tilde{X}\times_{\pp^1}\tilde{X}$.
We will thereby create a splitting of the Leray filtration, and if
$n$ is even or $\mbox{char }k\neq 2$ a section of the ring
epimorphism $\mbox{res}_{\cal Y}$.

\begin{lemma}\label{proyerelat} Assume $C({\cal Y})$.  Let
$\pi'^i\in A^n(\tilde{X}\times \tilde{X})$ be liftings of
$\pi^i_{\cal Y}$. Then $\pi'^i$ are such that
$$\pi'^i|Gr_F^\epsilon H^j(\tilde{X})=\delta_{i,j-\epsilon}$$ for
all $0\leq i, j\leq 2n-2$ and $\epsilon=0,1,2.$  
 The restriction of $\pi'^j$ to $F^2H^*(\tilde{X})$ is a
 projector which
 yields $0$ on $F^2H^j(\tilde{X})$ if $j\neq i+2$ and the identity if $j=i+2$.  Thus the
 restriction to $F^2_\rho$ is clearly independent of the lifting chosen.
\end{lemma} {\bf Proof:} The proof is laid out in \ref{accionss}.  
If $\pi'^i\in {\cal A}$ restricts to $\pi^i_{\cal Y}$, then
$\pi'^i|R^k\rho_*\Q_\ell=\delta_{i,k}$ for all $i,k$.
 Applying $H^\epsilon(\pp^1,\bullet)$ the Lemma follows.$\blacksquare$

    Whatever the choice of liftings $\pi'^i$, these operators commute with the
  K\"unneth projectors of $\tilde{X}$ by Lemma \ref{adjuntoh};
  this justifies the following definitions, which make sense under condition {\bf (A)} of Section
  \ref{lapiceslefschetz}. We define $\pi^{i,0}$ after the relation
  \newline
  $\pi^i_{\cal Y}=\,^t\pi_{\cal Y}^{2n-2-i}.$
\begin{caja}{Notation-Definition.}  Let $\pi^{i,2}$ denote the orthogonal projection onto $F^2H^{i+2}(\tilde{X})$,
and $\pi^{i,0}$ denote the transpose $^t\pi^{2n-2-i,2}$.  We define $\pi^{n-1,1}:=\pi_{\tilde{X}}^n-\pi^{n,0}-\pi^{n-2,2}$ and $\pi^{i,1}=0$ otherwise.  Then the  $\pi^{i,\epsilon}$ form a complete orthogonal system of projectors,
 and provide a splitting for the Leray filtration $F_\rho^\bullet$ of $H^*(\tilde{X}).$
\end{caja}

  The following Proposition is the relative equivalent to Lemma
  \ref{cx}.
\begin{prop}\label{hacherho} Let $H_\rho=\sum (n-1-i)\pi^i_\rho.$  Then $H_\rho$ is (characterised as) the only semisimple (algebraic) skew-symmetric operator supported on $\tilde{X}\times_{\pp^1}\tilde{X}$ mapping to $H_{\cal Y}$ under the specialisation map.  The
complete orthogonal system of projectors
$\{\pi^i_\rho\pi^j_{\tilde{X}}\}$ yields a splitting of the Leray
filtration $F_\rho^\bullet H^*(\tilde{X}).$ \end{prop} {\bf
Proof:} The proof is elementary. Let us prove existence first. The
correspondence $H_{\cal Y}=\sum (n-1-i)\pi^i_{\cal Y}$ lifts to a
non-unique correspondence $\tilde{H}$ on
$\tilde{X}\times\tilde{X}$ supported on
$\tilde{X}\times_{\pp^1}\tilde{X}$, which we may assume
skew-symmetric.  Indeed, $^tH_{\cal Y}=-H_{\cal Y}$ implies that
$^t\tilde{H}+\tilde{H}\in {\cal J}$ is nilpotent of order $3$ by
Proposition \ref{generico}.  Now, the minimal polynomial of
$\tilde{H}$ divides $R(x)=P(x)^3$, where $P(x)=x\prod_{i=1}^{n-1}
(x^2-i^2).$ The fact that $R(x)$ is odd implies the following.
\vskip.2cm {\bf Claim.} The semisimple part of the Jordan
decomposition of $\tilde{H}$ is skew-symmetric.
\begin{caja}{Proof of the Claim:} Write $$R_i(x)=\prod_{k\in [-n+1, n-1], k\neq i} (x-k)^3.$$  Write $$1=\sum_{-n+1}^{n-1} R_i(x)a_i(x),$$ with $a_i(x)$ quadratic polynomials.  It is clear that $R_i(-x)a_i(-x)=R_{-i}(x)a_{-i}(x).$  Multiplying the above by $x$, one has $x=\sum i R_i(x) a_i(x)+\sum R_i(x) a_i(x)(x-i)$, which yields the Jordan decomposition of $\tilde{H}$.  Substituting $\tilde{H}$ into $x$, it follows that the first sum is skew-symmetric and semisimple (and, of course, algebraic).
\end{caja}

 Let $H'_\rho$ be a semisimple, algebraic, skew-symmetric lifting of $H_{\cal Y}$. By Proposition \ref{generico}, $H'_\rho$ agrees with $H_\rho$ on $F^2H^*(\tilde{X})=\mbox{Im }\sum \pi^{i,2}.$   Transposing yields $$(H'_\rho-H_\rho)|\mbox{Im }\sum \pi^{i,0}=0.$$  It remains only to check equality on $\mbox{Im }\pi^{n-1,1}$: now $H'_\rho, H_\rho$ are nilpotent on $\mbox{Im }\pi^{n-1,1}$, hence $0$ by semisimplicity, thus completing the proof.
$\blacksquare$

\begin{corollary}\label{pirelat} With notations and assumptions of Proposition \ref{hacherho}, the relative projectors $\pi^i_\rho$ are the projections onto the primary components of
the operator $H_\rho$, and $\pi^i_\rho= ^t\pi_\rho^{2n-2-i}$.
Moreover,
$$\pi^i_{\tilde{X}}=\pi^{i,0}+\pi^{i-1,1}+\pi^{i-2,2} \mbox{  and  }\pi^i_\rho=\pi^{i,0}+\pi^{i,1}+\pi^{i,2},$$ where $\pi^{i-1,1}=0$ for $i\neq n$.\end{corollary}
{\bf Proof:}  The first assertion follows from Proposition
\ref{hacherho}, and $\pi^i_\rho$ are
thus polynomials in $H_\rho$. 
 The second assertion follows from the fact that $H_\rho$ is semisimple
skew-symmetric.  The rest follows from the Leray spectral sequence
of $\rho$, condition {(A)} and Lemma \ref{proyerelat}.
$\blacksquare$

  \begin{caja}{Observation-Definition} Let $\tilde{u}$ be a correspondence of degree $r$ of ${\cal Y}.$ Then
  \be \label{correspypi}
\tilde{u}=\sum \pi^{i+2r}_{\cal Y} \tilde{u} \pi^i_{\cal Y}.
  \ee
    This goes along with (and in fact implies) the commutation relation in
    Lemma \ref{adjuntoh}.  We define for each $u\in {\cal A}$ the
    following element of ${\cal A}$:
    \be\label{liftingchoice}
u_\rho:= \sum \pi^{i+2r}_\rho u \pi^i_\rho.
    \ee
  It is clear by construction that $u_\rho-u\in {\cal J}.$ If $u$
is a correspondence of degree $r$ on ${\cal Y}$, we will define
$u_\rho$ to be $u'_\rho$ for $u'$ a lifting of $u$ in ${\cal
A}$.  Later we will see that this definition is consistent.
\end{caja}
\begin{lemma}\label{rhochoice} The map $\psi:{\cal A}\rightarrow {\cal A}$ defined by $u\mapsto u_\rho$ satisfies ${\cal J}=\mbox{Ker
}\psi;$ in other words, $u_\rho=v_\rho$ if and only if $u-v$
induces $0$ on $Gr^\bullet_{F_\rho}$.  The image $\psi({\cal A})$
in degree $r$ consists of the $w \in {\cal A}$ such that
$\pi^{i+2r}_\rho w=w \pi^i_\rho$ for all $i$.  The map $\psi$ is a
linear projector which induces a section $\sigma$ of the natural
quotient map ${\cal A}\rightarrow {\cal A}/{\cal J}$ \`a la
Wedderburn-Malcev, and commutes with
 transposition.  As a result we have a well-defined ring homomorphism
 $$A^{n-1+*}({\cal Y}\times{\cal Y}) \rightarrow \psi({\cal A})$$ defined by $u\mapsto
 u_\rho,$ which agrees with the homomorphism $\psi \circ proj_{\cal J}
 \circ \varphi.$
\end{lemma}
{\bf Proof:} It is clear that ${\cal J}=\mbox{Ker }\psi.$  The
image of $\psi$ is easily characterised as the subspace of $u$
such that $\psi(u)=u$ (easily seen to agree with the description
$u\pi^i_\rho=\pi^{i+2r}_\rho u$ if $u$ is of degree $r$), whence
$\psi^2=\psi.$  By Corollary \ref{pirelat},
$\psi(\,^tu)=\,^t\psi(u)$. Finally, the terms $v_\rho\circ u_\rho$
and $(v\circ u)_\rho$ differ by
  an element of ${\cal J}\cap \mbox{Im }\sigma_{\cal Y}=(0)$, thus proving that $\psi$ is a ring
  homomorphism.  $\psi$ clearly induces a section ${\cal A}/{\cal
  J} \rightarrow {\cal A}$ of the quotient map, which gives rise to the map $u\mapsto u_\rho$ with target $A^{n-1+*}({\cal Y}\times{\cal Y})$. $\blacksquare$

\section{A relative ${\mathfrak sl}_2$-triple}

  We have obtained a set of relative K\"unneth projectors under
  the hypothesis $C({\cal Y}).$  In this section we assume $B({\cal Y})$ and we construct
  relative operators $^c\Lambda_\rho, \Lambda_\rho$ lifting $\,^c\Lambda_{\cal Y}, \Lambda_{\cal Y}$; this will
  give rise to a relative ${\mathfrak sl}_2$-triple $^c\Lambda_\rho, L_\rho,
  H_\rho$ whose action on $H^*(\tilde{X})$ will be exploited later.

  \begin{prop}\label{lambdarelcero}The following assertions hold.
\begin{enumerate}
    \item[(1)]  For any lifting $^c\Lambda'$ of $^c\Lambda_{\cal Y}$, the correspondence $^c\Lambda_\rho=\sum
  \pi^{i-2}_\rho\,
^c\Lambda'\pi^i_\rho$ is symmetric and independent of the lifting
$^c\Lambda'$ chosen.
    \item[(2)]  The operator $^c\Lambda_\rho$ satisfies $^c\Lambda_\rho \pi^{i,2}\subset \mbox{Im }\pi^{i-2,2}$ and
$^c\Lambda_\rho \pi^{i,0}\subset \mbox{Im }\pi^{i-2,0}$.  In fact
$^c\Lambda_\rho \pi^{i,0}=\pi^{i-2,0}\,^c\Lambda_\rho$ and
$^c\Lambda_\rho \pi^{i,2}=\pi^{i-2,2} \,^c\Lambda_\rho.$
\end{enumerate}
  \end{prop}
{\bf Proof:} \textbf{(1)} follows from Lemma \ref{rhochoice}.
\textbf{(2)} follows directly from Proposition \ref{fil} and Lemma
\ref{rhochoice}. $\blacksquare$

\begin{lemma}\label{piezasola}  $\mbox{Im }\pi^{n-1,1}=P^n(X)\oplus
V(\Delta)(-1)$ (compare Katz ~\cite{SGA7} Exp. XVIII Th. 5.7) and
$\mbox{Im }\pi^{n,0}$ is the image of $\Delta(H^{n-2}(Y))$ via the
inclusion
$$\iota_*\oplus h^*:H^{n-2}(Y)(-1)\oplus H^{n-2}(Y)(-1) \hookrightarrow H^n(X)\oplus H^{n-2}(\Delta)$$
given by the decomposition of Proposition \ref{katzuno}.
 On the other hand, $$H^n(X)\cap F^1_\rho H^*(\tilde{X})=P^n(X)\oplus 0.$$\end{lemma}
{\bf Proof:} The sought-for image of $\pi^{n-1,1}$ coincides with
the orthogonal in $F^1_\rho H^n(\tilde{X})$ of
$k_*H^{n-2}(Y)=\mbox{Im }\pi^{n-2,2}.$

  Note the orthogonal decomposition
  \be
H^n(\tilde{X})=(P^n(X)\oplus V(\Delta))\oplus^\perp
(\iota_*H^{n-2}(Y) \oplus h^*H^{n-2}(Y)(-1)).
  \ee
  The piece $P^n(X)\oplus V(\Delta)$ is clearly within $F^1_\rho$
  and orthogonal to $F^2_\rho H^n(\tilde{X})$; by a dimension count (see Corollary \ref{dualidades}) we have $(P^n(X)\oplus V(\Delta))\oplus^\perp F^2_\rho H^n(\tilde{X})=F^1_\rho H^n(\tilde{X});$
  the equality $H^n(X)\cap F^1_\rho=P^n(X)\oplus 0$ is thus established.

    We define
  $W=H^{n-2}(Y)(-1)^{\oplus 2}$ and view it as a quadratic
  subspace of $H^n(\tilde{X})$ via $\left(\begin{array}{cc} \iota_* & 0 \\ 0 &  h^* \end{array}
  \right)$.  Write $W=W_1\oplus W_2$ ($\oplus$ not
  orthogonal), where $W_1=\mbox{Im } \left(\begin{array}{c} 1 \\ -1 \end{array}
  \right)$ represents $F^2_\rho  H^n(\tilde{X})$ and $W_2=\mbox{Im
  }\left(\begin{array}{c} 1 \\ 1 \end{array}
  \right)$ -- note that both $W_i$ are self-orthogonal.  We will show that $\pi^{n,0}$ is given by the projection onto $W_2$.

Let $w=(w_1, w_2), w'=(w'_1, w'_2) \in W=W_1\oplus W_2$.  Then
$\langle w_i, w'_i \rangle=0,$ and \bean \langle w_1,
w'_1+w'_2\rangle=\langle w_1+w_2, w'_2 \rangle,
 \eean
which shows $\,^t\pi^{n,0}=\pi^{n-2,2}.$
    The Lemma is thus established. $\blacksquare$

\begin{lemma}\label{primirelativos} Notations and assumptions as above, let $i\leq n-1.$  Then
$$\mbox{Im }\pi^{i,0}=(P^i(X)\oplus 0)\oplus (\iota_* \oplus h^*)H^{i-2}(Y).$$

  \end{lemma}
  {\bf Proof:} The first assertion is similar to Lemma
  \ref{piezasola}.  The piece $P^i(X)\oplus 0\subset \mbox{Im }\pi^{i,0},$ since the
  image of the projector $f^*\,^tp^i_Xf_*$ is contained in
  $\mbox{Im }\pi^{2n-i-2}$ by the equality $k_*L^{n-i-1}P^i(Y)=L^{n-i}P^i(X)\oplus 0;$ the piece
  $(\iota_*\oplus h^*)\Delta(H^{i-2}(Y))\subset \mbox{Im }\pi^{i,0}$ by a similar argument to Lemma \ref{piezasola}.  To
  prove the second assertion we note the following: if
  $x\in P^i(X)$ for $i\leq n-2,$ then $L^{n-i-1}P^i(X)\oplus 0
  \cap
   k_*H^{2n-4-i}(Y)=0,$ since \newline $j^*:P^i(X)\hookrightarrow P^i(\Delta)$ is
  injective; here we have used Lemma \ref{calculostildes2}.  The case $i=n-1$ is
  obvious.
   $\blacksquare$

  Finally we obtain the desired ${\mathfrak sl}_2$-triple.
\begin{prop}\label{sltriple}  We have a relative ${\mathfrak sl}_2$-triple $^c\Lambda_\rho, L_\rho, H_\rho.$  A relative Lefschetz
 isomorphism holds:
 $$L^i_\rho:\mbox{Im }\pi^{n-1-i}_\rho \rightarrow \mbox{Im }\pi^{n-1+i}_\rho$$ for $1\leq i\leq n-1$. The projectors $p^{i}_\rho$ are algebraic for $i\leq n-1$,
 and we have symmetric operators $p^{n-1+j}_\rho$ derived from $p^{n-1+j}_{\cal
 Y}$ for $0\leq j\leq n-1$.  The map $u\mapsto u_\rho$ yields an isomorphism of rings $\Q\langle L_\rho, \Lambda_\rho\rangle\cong \Q\langle L_{\cal Y}, \Lambda_{\cal
 Y}\rangle$ which preserves transposition.
 \end{prop}
{\bf Proof:}  The ${\mathfrak sl}_2-$identities $[H_\rho,
^c\Lambda_\rho]=2 ^c\Lambda_\rho$, $[H_\rho,L_\rho]=-2L_\rho$ and
$[^c\Lambda_\rho, L_\rho]=H_\rho$ and the isomorphism between
$\Q\langle L_\rho, \Lambda_\rho\rangle$ and $\Q\langle L_{\cal Y},
\Lambda_{\cal
 Y}\rangle$
follow from Lemma \ref{rhochoice} and Corollary \ref{rhocardos}.
The operators $\frac{1}{m^i}L_{\tilde{X}}^i, L^i$ and $L_\rho^i$
induce the same map on $\mbox{Gr}_{F_\rho}^\bullet$ by Proposition
\ref{leray} and Lemma \ref{calculostildes2}. The `relative
Lefschetz isomorphism' can be checked by passing to
$\mbox{Gr}_{F_\rho}$, or simply by using the identities
$(\Lambda_\rho^{i}L^{i}_\rho-1)\pi^{n-1-i}_\rho=0$ for $i<n-1.$
 $\blacksquare$

  We now can view the Lefschetz theory of ${\cal Y}$ within
  $H^*(\tilde{X}).$

\begin{prop}\label{primirelatpeques} The relative primitive
projectors $p^i_\rho$ for $i\leq n-2$ are described as follows.
$p^i_\rho|H^i({\tilde{X}}=f_* p^i_X f^*,$ $p^i_\rho
H^{i+1}(\tilde{X})=0$ and $\mbox{Im }p^{i+2}_\rho
H^{i+2}(\tilde{X})=k_*P^i(Y)=LP^i(X)\oplus P^{i-2}(\Delta)\cap
F^2_\rho.$
\end{prop}
{\bf Proof:}  By Lemma \ref{primirelativos} and Corollary
\ref{inicioprim}, $P^i(X)\oplus 0\subset \mbox{Im }\pi^{i,0}$ is
isomorphic to $H^0({\cal P}^i) \subset H^0({\cal R}^i)$ via the
obvious map.  Since ${\cal P}^i$ is constant, the image of
$p^i_\rho$ in degree $i+1$ is zero; the last assertion follows
from Corollary \ref{inicioprim}. $\blacksquare$
\begin{prop}\label{penemenos} 
  The projector $p^{n-1}_\rho$ satisfies the following properties:
\begin{enumerate}
    \item $\pi^{n-1}_\rho p^{n-1}_\rho=p^{n-1}_\rho$ and $p^{n-1}_\rho= \,^tp^{n-1}_\rho$;
    \item $\pi^{n-1,1} p^{n-1}_\rho= \pi^{n-1,1}$, i.e. the
    orthogonal projection onto $$P^n(X)\oplus V(\Delta)\subset
    P^n(\tilde{X});$$
    \item $\pi^{n-1,2}p^{n-1}_\rho$ is the orthogonal projection onto
    $LP^{n-1}(X)\oplus 0$, and $\pi^{n-1,0}p^{n-1}_\rho$ is the
    projection onto $P^{n-1}(X)\oplus 0$.
\end{enumerate}
In all, 
the projector $p^{n-1}_\rho$ can be expressed as

\be \label{formulapenemenos} p^{n-1}_\rho=f^*p^{n-1}_X
f_*+\pi^{n-1,1}+f^*\,^tp^{n-1}_X f_*, \ee and $f_* p^{n-1}_\rho
f^*=p^{n-1}_X+p^n_X+\,^tp^{n-1}_X.$
\end{prop}
{\bf Proof of Proposition \ref{penemenos}:}  
\begin{enumerate}
    \item is a straightforward consequence of Proposition \ref{sltriple}.
    \item Since $R^{n-3}\rho_*\Q_\ell$ is constant, one has
    $H^1(R^{n-1}\rho_*\Q_\ell)=H^1({\cal P}^{n-1}_\rho)$; the rest
    follows from Lemma \ref{piezasola}.
    \item The computation $H^2({\cal P}^{n-1}_\rho)=LP^{n-1}(X)\oplus 0$ follows from
    Corollaries \ref{lefschetztildedec},
    \ref{inicioprim}.

      Thus $\mbox{Im }\pi^{n-1,2}p^{n-1}_\rho=LP^{n-1}(X)\oplus
      0.$ Using the Poincar\'e duality pairing yields $\mbox{Im }\pi^{n-1,0}p^{n-1}_\rho=P^{n-1}(X)\oplus
      0$, by  Lemma \ref{primis} and Corollary \ref{lefschetztildedec}.
      $\blacksquare$
\end{enumerate}

\section{The Main Theorem}

  This section is devoted to proving the following result.
\begin{caja}{Main Theorem.} Let $X$ be a smooth, projective variety of
dimension $n$. Assume the Lefschetz standard conjecture for the
generic fibre ${\cal Y}/k(t)$ of a Lefschetz pencil satisfying
{\bf (A)}.  Then $\Lambda-p^{n+1}$ is algebraic.
\end{caja}

  We will prove this result in a series of steps, obtaining the algebraicity of the K\"unneth projectors $\pi^i_{\tilde{X}}$ for $i\neq n-1,n, n-1$ in the course of our proof.

\subsection{The algebraicity of some projectors}
    We start by proving the following.
\begin{prop}\label{bandc}
$B({\cal Y})$, with ${\cal Y}$ as above, implies the algebraicity
of the K\" unneth projectors $\pi^i_X$ for $i\leq n-2$ (hence that
of $\pi^i_X$ for $i\geq n+2$) and that of the primitive projectors
$p^0,\ldots, p^{n-2}$.
\end{prop}
A couple of lemmas will be required to establish this Proposition.
\begin{lemma} \label{lemabc} The following statements hold.
\begin{enumerate}
\item[(i)  ] The identity $\iota_*H^{i-2}(Y)=LH^{i-2}(X)$ holds
for all $0\leq i\leq 2n$, and \linebreak
$\iota_*:H^{i-2}(Y)\rightarrow H^i(X)$ is injective for $i\leq
 n.$  For all $i\leq n$, $$\mbox{Im}(\pi^i_X-p^i_X)=LH^{i-2}(X).$$
\item[(ii) ] For all $i>n$,
$\mbox{Im}(\pi^i_X-\,^tp^{2n-i}_X)=L^{i-n+1}H^{2n-i-2}(X).$
 \item[(iii)] Suppose $B(Y)$ holds for $Y$ a smooth hyperplane section of $X$.  Then for $0\leq i\leq n$,
$$\pi^i_X-p^i_X=\iota_*\Lambda_Y\pi_Y^i\iota^*$$ is algebraic.  Thus the transposed operators $\pi^{2n-i}-L^{n-i}p^{2n-i}$ are algebraic for $0\leq i\leq n$.
\item[(iv) ]  The hypothesis $B({\cal Y})$ of the Main Theorem
implies $B(Y)$ for a suitable hyperplane section.
\end{enumerate}
\end{lemma}
{\bf Proof:} Statement \textit{(iv)} follows by Proposition
\ref{basiclambda}\textit{(2)} and specialisation.  The rest is
straightforward. $\blacksquare$

 We consider a suitable Lefschetz pencil for $X$, and prove the algebraicity of $\pi^i_X$ for $i\leq n-2$.

   It suffices to prove
  that the operators $\pi_{\tilde{X}}^i$ are algebraic for $i=0, \cdots, n-2$, since $\pi^{2n-i}=\,^t\pi^i$ (Kleiman ~\cite{Kleiman} showed already that $\pi^0, \pi^1$ are algebraic in general).
\begin{lemma}\label{ctrivial} The projectors $\pi^{i-2,2}$ are algebraic for $i\leq n$, and so are $\pi^{i-2,0}$.
\end{lemma}
{\bf Proof of Lemma \ref{ctrivial}:}
\begin{itemize}
\item  We know that $k_*=\iota_* \oplus -h^*$.  Let $i\leq n$.
Let us prove that
  $$\pi^{i-2,2}=\pi^{i-2,2}f^*(\pi^i_X-p_X^i)f_*\pi^{i-2,2}.$$  Indeed, by
  Lemma \ref{lemabc}, the image of $k_*(y)=\iota_*(y)\oplus
  -h^*(y)$ by $(\pi^i_X-p_X^i)\pi^{i-2,2}$ is $\iota_*(y)\oplus
  0$ if $y\in H^{i-2}(Y)$ and $0$ otherwise.  Applying $\pi^{i-2,2}$
  to $\iota_*(y)\oplus 0$ yields $\iota_*(y)\oplus -h^*(y)$ by
  Lemma \ref{lemabc}\textit{(i)}.
  \item Rewriting the previous step we get $$\pi^{i-2,2}=\pi^{i-2}_\rho f^*(\pi^i_X-p_X^i)f_*\pi^{i-2}_\rho.$$
\item By the above, $\pi^{i-2,2}$ is algebraic for $i\leq n-2$,
and the operator $\pi^{i-2}_\rho-\pi^{i-2,2}=\pi^{i-2,0}$ is
algebraic for $i\leq
n-2.$ The Lemma is thus settled. $\blacksquare$ 
\end{itemize}
\begin{lemma}The projectors $\pi^i_{\tilde{X}}$ are algebraic for $i\neq n-1, n, n+1$. \end{lemma}
{\bf Proof:} The proof is immediate, since for $i\leq n-2$ the
operator $\pi^i_{\tilde{X}}=\pi^{i,0}+\pi^{i-2,2}$ is algebraic by
Lemma \ref{ctrivial}. $\blacksquare$

{\bf Proof of Proposition \ref{bandc}:}
   It remains only to check that  $p^i_X=\pi^i_X-(\pi^i_X-p^i_X)$ is algebraic for
$i\leq n-2$; this holds by Lemma \ref{lemabc}\textit{(iii)}.
Proposition \ref{bandc} is thus established. $\blacksquare$

\subsection{Proof of the Main Theorem}\label{subsectmainth}

  We finally prove the Main Theorem.

  Assume that $B({\cal Y})$ holds for ${\cal Y}$ the generic fibre
  of a Lefschetz fibration $\rho$ of $X$ satisfying condition {\bf (A)}.

  By Lemma \ref{prinduccion}, we have the following identity:
  \be\label{casifinal}
\iota^*(\Lambda_X-p^{n+1}_X)-\Lambda_Y\iota^*= \sum_{j=n+2}^{2n-2}
\iota^*L^{j-n-1} p^j_X.
 \ee
\begin{caja}{Aside (not necessary):}  Assuming $B(Y)$, the l.h.s. of is algebraic if and only if $\Lambda_X-p^{n+1}_X$ is.  This follows from the identity
$^t[\iota^*(\Lambda_X-p^{n+1}_X)]
\iota^*(\Lambda_X-p^{n+1}_X)=\Lambda_X-p^{n+1}_X.$
\end{caja}
   The next step is to prove that the r.h.s. of (\ref{casifinal}) is
  algebraic.  This will follow from the next Lemma ($j=2n-i$).

\begin{lemma}\label{lemafinal} Assume $B({\cal Y})$.  For $0\leq i\leq n-2$,
the operator
$$L^{n-i-1}p^{2n-i}_X=\Lambda_X\,^tp^{i}_X=f_*\Lambda_\rho f^*\,^tp_X^{i}$$ is algebraic.\end{lemma}
{\bf Proof of Lemma \ref{lemafinal}:} Let $i\leq n-2$; then
$p^i_X$ is algebraic by Proposition \ref{bandc}.
 Consider the subspace $W=(L^{n-1-i}P^i(\tilde{X})\oplus L^{n-1-i}P^{i-2}(\Delta))\cap F^2_\rho$,
which agrees with the image of
$$k_*=\iota_*\oplus (-h^*):L^{n-2-i}P^i(Y)\rightarrow L^{n-1-i}P^i(X)\oplus
L^{n-2-i}P^i(\Delta).$$  The first component is an isomorphism,
and the second is injective, being bijective if $i<n-2.$  On
applying $L$, which coincides with $L_\rho$ on $F^2_\rho$, we have
an isomorphism $$L:W \overset{\sim}{\rightarrow}
L^{n-i}P^i(X)\oplus 0 \subset F^2_\rho;$$ the piece
$L^{n-i}P^i(X)\oplus 0=L^{n-i}_{\tilde{X}}P^i(\tilde{X})$ equals
$k_*L^{n-1-i}_Y P^i(Y)$ -- see Corollary \ref{lefschetztildedec}.
$L$ is thus an isomorphism between $W'$ and $L^{n-i}P^i(X)\oplus
0$ (by Corollary \ref{lefschetztildedec}, $m\cdotp L$ and
$L_{\tilde{X}}$ agree on $F^2_\rho$). The identity \bean
 L_{\cal Y}\Lambda_{\cal Y}= 1_{\cal Y}-\sum_{i=0}^{n-1} p^i_{\cal Y} \eean
 (~\cite{Kleiman} p. 372) translates by Proposition \ref{sltriple} into
$$L_\rho\Lambda_\rho=1_{\tilde{X}}-\sum_{i=0}^{n-1} p^i_{\rho},$$
thus showing that $\Lambda_\rho$ defines the inverse isomorphism
to $L:W \rightarrow L^{n-i}P^i(X)\oplus 0.$ Taking the
$X$-component yields the inverse
$$L^{n-i}P^i(X)\oplus 0 \rightarrow W
\rightarrow L^{n-1-i}P^i(X)$$ of $L$, which coincides with
$\Lambda_X|L^{n-i}P^i(X)$ -- here we have used that $L$ agrees
with $L_\rho$ on $F^2_\rho$, and that $p^i_\rho$ acts as $0$ on
$H^j(\tilde{X})$ for $j\geq n+2$. We have thus proven that
$\Lambda_X\,^tp^{n-i}_X=L^{n-1-i}p^{n+i}=f_*\Lambda_\rho
 f^*\,^tp^{n-i}_X$ is algebraic by Propositions \ref{sltriple} and
\ref{bandc}. $\blacksquare$

\begin{lemma} \label{finalsi} Assuming the hypotheses of the Main Theorem, the operator $$\Lambda_X\pi^i_X=\pi^{i-2}_X\Lambda_X=\Lambda_X (\pi^i_X-p^i_X)$$
 is algebraic for $i\leq n.$
  The operator $(\pi^{n-1}_X-p^{n-1}_X)\Lambda_X=(\Lambda_X-p_X^{n+1})\pi^{n+1}_X$ is algebraic.\end{lemma}
{\bf Proof:} The identities are clear; let us prove algebraicity
of the above operators.  By Lemmas \ref{prinduccion} and
\ref{lemabc}\textit{(iii)} we have
$$\Lambda_X
(\pi^i_X-p^i_X)=\Lambda_X\iota_*\Lambda_Y\pi^i_Y\iota^*=\iota_*\Lambda_Y^2\pi^i_Y\iota^*
+\left(p_X^{n+1}+\sum_{j=n+2}^{2n-2}
p_X^{j}L^{j-n+1}\right)\iota_*\Lambda_Y\pi^i_Y\iota^*=$$
$$=\iota_*\Lambda_Y^2\pi^i_Y\iota^*
+\left(\sum_{j=n+2}^{2n-2}
p_X^{j}L^{j-n+1}\right)\iota_*\Lambda_Y\pi^i_Y\iota^*,$$ which is
algebraic for $i\leq n+1$ by Lemma \ref{lemafinal}, thereby
establishing the Lemma. $\blacksquare$

\begin{caja}{Proof of Main Theorem:} Under the hypotheses of this Section, the operator $\Lambda_X-p^{n+1}_X$ is algebraic.\end{caja}
Indeed, we have proven in Lemma \ref{finalsi} that
$\pi^{n-2}_X\Lambda_X=\Lambda_X\pi^n_X$ and
$\pi^{n-3}_X\Lambda_X=\Lambda_X\pi^{n-1}_X$ are algebraic, as well
as the algebraicity of
$\left(\Lambda_X-p^{n+1}_X\right)\pi^{n+1}_X$.  It now remains to
establish the algebraicity of $\Lambda_X
\sum_{i=n+2}^{2n}\pi_X^i$. Again, Lemma \ref{finalsi} shows that,
for $r\geq 2$, the operator
$\,^t\left(\Lambda_X\pi_X^{n+r}\right)=\pi_X^{n-r}\Lambda_X$ is
algebraic.  On the other hand,
$(\Lambda_X-p^{n+1}_X)\pi_X^{n+1}=\Lambda(\pi_X^{n+1}-\,^tp_X^{n-1})$.
Altogether this shows that the operator
$$\Lambda_X-p^{n+1}_X=\Lambda_X\left(\sum_{k=0}^{n}\pi^i_X+(\pi_X^{n+1}-\,^tp^{n-1}_X)+\sum_{k=n+2}^{2n}\pi^i_X\right)$$
is algebraic.  The Main Theorem is thus established.
$\blacksquare$
\subsection{Final comments} \label{finalcomments}

  On the field of definition of the correspondences $\,^c\Lambda,
  \Lambda, \pi^i$ we would like to say the following.  The correspondence $L$ is $k$-defined.  Now assume that $k$ is
  perfect: the operator $H$ and all K\"unneth projectors are Galois invariants; if they are
  algebraic, an elementary argument furnishes $k$-defined
  algebraic representatives for $H,\pi^i$.  The operator
  $\,^c\Lambda$ is uniquely determined by its ${\mathfrak sl}_2$-partners $H, L$ (see e.g. ~\cite{Andre04} 5.2.2), hence Galois invariant, and again is represented
  by a $k$-defined algebraic class if it is algebraic.  If $k$ is
  arbitrary, one needs to descend from a purely inseparable finite
  extension $k'$ to $k$; which is standard, since the natural map $X\times_k k' \rightarrow
  X$ is a homeomorphism.  Thus our form of the conjectures is in fact
  equivalent to that of ~\cite{G} ~\cite{Kleiman}, where it was
  assumed that $k=\overline{k}.$  This formulation was required in
  order to obtain algebraic cycles supported on $\tilde{X}\times_{\pp^1}\tilde{X}.$

\begin{caja}{Restatement of the Main Theorem:} In the language of Proposition \ref{basiclambda},
our Main Theorem shows precisely that $\theta^i$ is induced by an
algebraic cycle (to wit $(\Lambda-p^{n+1})^{n-i}$)  for $i\leq
n-2$. By the proof of ~\cite{Kleiman} Lemma 2.4, one derives that
$\pi^i_X$ is algebraic for $i\neq n-1, n, n+1.$
\end{caja}

\begin{caja}{On the results needed in the Proof:} The path we have travelled in order to settle our Main Theorem is not the shortest
possible.
 Lemma \ref{lemafinal} does not need $\Lambda_\rho$ or
$\pi^i_\rho$, but
 merely arbitrary liftings of $\Lambda_{\cal Y}, \pi^i_{\cal Y}$ to $\tilde{X}\times_{\pp^1}\tilde{X}$, since the proof of Lemma \ref{lemafinal}
 requires only working on $F^2_\rho$.  Lemma \ref{finalsi} relies
 solely on Lemma \ref{lemafinal} and material from Section \ref{seccgeneral}, and the Main Theorem rests on Lemma \ref{finalsi}.
     The use of Lemma \ref{finalsi} allows for a direct
proof of the algebraicity of $\Lambda-p^{n+1}$ without Proposition
\ref{bandc}, but then one should use $\theta^0, \cdots ,
\theta^{n-2}$ as in the proof of ~\cite{Kleiman} Lemma 2.4, to
derive that $\pi^i_X$ is algebraic for $i\neq n-1, n,
n+1$.\end{caja}
\begin{caja}{On the support of $\pi^i_{\tilde{X}}$:} We have proven that $\pi^{i-2,2}$ is algebraic for $i\leq n.$  However, our proof does not necessarily
imply that $\pi^{i-2,2}$ is supported on
$D=\tilde{X}\times_{\pp^1}\tilde{X}.$  In fact, it follows from
Remark \ref{remarcagenerico} and the proof of Proposition
\ref{generico} that $\pi^{i-2,2}$ is not supported on $D$ unless
it is $0$ (\textit{i.e.} when $H^{i-2}(X)=0$).\end{caja}
\subsection{The operator $p^{n+1}_X$}

  This section is a complement to the Main Theorem.
 The operator $p^{n+1}_X$ (and so $\iota^*p_X^{n+1}$) is of central importance in the Lefschetz theory of $X$.
\begin{lemma}  Assume $B({\cal Y})$ for ${\cal Y}$ the general fibre of a Lefschetz pencil of $X$ satisfying {\bf (A)}.  The algebraicity of $p^{n-1}_X$ implies that of $p^n_X$ and the conjecture $C(X)$.
\end{lemma} The lemma follows from Proposition \ref{penemenos}. $\blacksquare$

\begin{lemma}\label{pnmas1} Let $X$ be a projective smooth, $n$-dimensional variety.  The operator $\iota^*p^{n+1}$ is algebraic if $p^{n+1}$ is. \end{lemma}
{\bf Proof:} The result stems from the following identity: \be
\label{penemas1} \,^t(\iota^*p^{n+1})
\iota^*p^{n+1}=p^{n+1}Lp^{n+1}=p^{n+1}.\blacksquare \ee

  Enclosed in $p^{n+1}$ is information about the space of
  vanishing cycles, and also $p^{n-1}$.  For instance, as we shall
  see below,
  the algebraicity of $p^{n+1}$ allows us to speak of the
  \textbf{motive of vanishing cycles}.

\begin{prop} With the above notations, the following statements
hold.
\begin{enumerate}
    \item The algebraicity of $p^{n+1}_X$ implies that of $p^{n-1}_X$.
    \item Let $V(Y)$ be the space of vanishing cycles of a smooth hyperplane section
    $Y$, and $e_{V(Y)}$ be the orthogonal projection $H^*(Y)\twoheadrightarrow V(Y) \hookrightarrow H^*(Y).$  Then $$p^{n-1}_Y -e_{V(Y)}=\iota^* p^{n+1}
    \iota_*.$$
    \item If $B(Y)$ holds and $p^{n+1}$ is algebraic, then so is $e_{V(Y)}$.
\end{enumerate}$\blacksquare$
\end{prop}

\begin{prop} The operator $p^{n+1}$ cannot be obtained as $f_*uf^*$, with $u$ an algebraic cycle supported on $D=\tilde{X}\times_{\pp^1}\tilde{X}.$  \end{prop}
{\bf Proof:} It follows from Proposition \ref{penemenos} that
$LP^{n-1}(X)\oplus 0\subset F^2_\rho$ and $$P^{n-1}(X)\oplus
H^{n-3}(\Delta) \cap F^1_\rho=0;$$ since every correspondence
supported on $D$ preserves the Leray filtration by Proposition
\ref{fil}, the assertion is clear. $\blacksquare$

{\bf Email address:} jjramon@maths.ucd.ie

\begin{thebibliography}{xxx}

\bibitem{Andre} \rm{Y. Andr\'e}, \textit{Pour une th\'eorie inconditionelle des motifs,} Publ. Math. IHES 83 (1996), 5-49.
\bibitem{Andre04}\rm{Y. Andr\'e,} \textit{Une introduction aux motifs,} Soc. Math. France, Paris 2004.

\bibitem{Artin} M. Artin, \textit{Grothendieck Topologies},
Harvard Univ. Seminar, 1962.
\bibitem{dJ} A.J. de Jong, \textit{Smoothness, semi-stability and alterations.}  Publ. Math. IHES (1996)
\bibitem{SGA7} P. Deligne, N. Katz (eds.), \textit{S\'eminaire de G\'eometrie Alg\'ebrique 7,} Vol. II.  Springer Lecture Notes in Mathematics 340.  Springer Verlag, Heidelberg, 1973.
\bibitem{DeWeI} P. Deligne, \textit{La Conjecture de Weil I},
Publ. Math. IHES 43 (1974), 273-307.
\bibitem{DeWeII} P. Deligne, \textit{La Conjecture de Weil II}, Publ. Math. IHES 52(1980), 313-428.
\bibitem{DMOS} \rm{P. Deligne, J.S. Milne, A. Ogus, K.Y. Shih}, \textit{Hodge cycles, motives and Shimura varieties}, LNM 900.  Springer-Verlag, Berlin-Heidelberg-New York, 1982.
\bibitem{FKiehl} E. Freitag, R. Kiehl, \emph{Etale cohomology and the Weil
conjecture}, Ergeb. der Math. Wiss., Springer-Verlag, Heidelberg
1988.
\bibitem{F} W. Fulton, \textit{Intersection Theory}, Springer-Verlag, Heidelberg 1984.

\bibitem{SGA4} \rm{A. Grothendieck} \textit{et al.}, \textit{S\'eminaire de G\'eometrie Alg\'ebrique 4},  LNM 269, 270,305. Springer-Verlag, Heidelberg 1972, 1973.
\bibitem{SGA5} \rm{A. Grothendieck} \textit{et. al.}, \textit{S\'eminaire de G\'eometrie Alg\'ebrique 5, 1965-66}. LNM
589. Springer-Verlag, Heidelberg 1977.
\bibitem{G} A. Grothendieck, \textit{Standard conjectures on Algebraic Cycles.} In: Bombay Colloquium, Oxford, 1969.
\bibitem{HAG} \rm{R. Hartshorne}, \textit{Algebraic Geometry},
GTM , Springer-Verlag, New York 1978.

\bibitem{Jannsen} \rm{U. Jannsen}, \textit{Motives, numerical equivalence, and semisimplicity.}  Inv. Math. 107 (1992) 447-452.
\bibitem{Seattle}  U. Jannsen, S. Kleiman, J. P. Serre (eds.) \emph{Motives}. Proc. Symp. Pure Math. 55 (Seattle, WA 1991), Part
1. AMS, Providence, RI 1995.
\bibitem{KM} N.M. Katz, W. Messing, \textit{Some consequences of the Riemann Hypothesis for varieties over finite fields.}  Inv. Math. 23
(1974), 73-77.
\bibitem{Kleiman} S. Kleiman, \textit{Algebraic cycles and the Weil conjectures.}  In: Dix Expos\'es sur la Cohomologie des Sch\'emas. A.Grothendieck, N. H. Kuiper (eds.), North-Holland, Amsterdam 1970.
\bibitem{Kleiman2} S. Kleiman, \textit{The Standard Conjectures}.  In ~\cite{Seattle}
3-20.
\bibitem{KleimanOslo} S. Kleiman, \textit{Motives}.  In: \textit{Algebraic Geometry,
Oslo 1970}, F. Oort (ed.), Walters-Noordhof, Groningen 1972, pp.
53-82.
\bibitem{Looijenga} E. Looijenga, \emph{Intersection Cohomology.}
  In: \emph{Complex Algebraic Geometry}, J. Koll\'ar, (ed.), AMS 1997.
\bibitem{Roberts} \rm{J. Roberts}, \textit{Chow's moving lemma}.  In: Algebraic Geometry, Oslo 1970.  F. Oort, ed.  Wolters-Noordhoff,  Groningen, 1972.
\bibitem{Scholl} \rm{A.J. Scholl}, \textit{Classical Motives}. In: ~\cite{Seattle} pp. 163-187.
\end{thebibliography}
\end{document}